\documentclass[12pt]{article}

\usepackage{amsmath, amssymb, booktabs, color, algorithm, algorithmic}
\usepackage{algorithm}

\usepackage{listings,color,url}

\definecolor{string}{rgb}{0,0.5,0.25}
\definecolor{comment}{rgb}{0,0.5,0}
\newtheorem{thm}{Theorem}[section]
\newtheorem{prop}{Proposition}[section]
\newtheorem{lem}{Lemma}[section]

\newtheorem{prob}{Problem}[section]

\newtheorem{assum}{Assumption}[section]

\newcommand{\argmin}{\operatornamewithlimits{argmin}}

\numberwithin{equation}{section}

\begin{document}
\makeatletter

\begin{center}
\large{\bf Fixed Point Algorithm for Solving Nonmonotone Variational Inequalities
in Nonnegative Matrix Factorization}\\
\small{This work was supported by the Japan Society for the Promotion of Science through a Grant-in-Aid for
Scientific Research (C) (15K04763).}
\end{center}\vspace{2mm}

\begin{center}
\textsc{Hideaki Iiduka} and \textsc{Shizuka Nishino}\\
Department of Computer Science, 
Meiji University,
1-1-1 Higashimita, Tama-ku, Kawasaki-shi, Kanagawa 214-8571 Japan\\ 
(iiduka@cs.meiji.ac.jp, shizuka@cs.meiji.ac.jp)
\end{center}

\vspace{2mm}

\footnotesize{
\noindent\begin{minipage}{14cm}
{\bf Abstract:}
Nonnegative matrix factorization (NMF), which is the approximation of a data matrix as the product of two nonnegative matrices, is a key issue in machine learning and data analysis. One approach to NMF is to formulate the problem as a nonconvex optimization problem of minimizing the distance between a data matrix and the product of two nonnegative matrices with nonnegativity constraints and then solve the problem using an iterative algorithm. The algorithms commonly used are the multiplicative update algorithm and the alternating least-squares algorithm. Although both algorithms converge quickly, they may not converge to a stationary point to the problem that is equal to the solution to a nonmonotone variational inequality for the gradient of the distance function.
This paper presents an iterative algorithm for solving the problem that is based on the Krasnosel'ski\u\i-Mann fixed point algorithm. Convergence analysis showed that, under certain assumptions, any accumulation point of the sequence generated by the proposed algorithm belongs to the solution set of the variational inequality. Application of the {\tt 'mult'} and {\tt'als'} algorithms in MATLAB and the proposed algorithm to various NMF problems showed that the proposed algorithm had fast convergence and was effective.
\end{minipage}
 \\[5mm]


\hbox to14cm{\hrulefill}\par

\section{Introduction}\label{sec:1}
Nonnegative matrix factorization (NMF) \cite{lee1999,paat1994} has attracted a great deal of attention since it has many potential applications, such as machine learning \cite{lee1999}, data analysis \cite{bru2004,cic2009,pau2004}, and image analysis \cite{gil2015}, and its complexity has been studied \cite{vava2009}. NMF is used to factorize a given nonnegative matrix $V \in \mathbb{R}_+^{M \times N}$ into two nonnegative factors $W \in \mathbb{R}_+^{M \times R}$ and $H \in \mathbb{R}_+^{R \times N}$ such that
\begin{align}\label{NMF}
V \approx WH
\end{align}
in the sense of a certain norm of $\mathbb{R}^{M \times N}$, where $R < \min \{M,N\}$ is a positive integer given in advance. The problem of finding $W \in \mathbb{R}_+^{M \times R}$ and $H \in \mathbb{R}_+^{R \times N}$ satisfying \eqref{NMF} can be formulated as a nonconvex optimization problem: 
\begin{align}\label{prob:0}
\text{Minimize } f (W,H) := \frac{1}{2} \left\| V - WH \right\|^2_{\mathrm{F}}
\text{ subject to } W \in \mathbb{R}_+^{M \times R}
\text{ and } H \in \mathbb{R}_+^{R \times N},
\end{align}
where $\| \cdot \|_{\mathrm{F}}$ stands for the Frobenius norm of $\mathbb{R}^{M \times N}$.

An algorithm commonly used for solving Problem \eqref{prob:0} is the {\em multiplicative update algorithm} \cite{lee1999,lee2001} (see Algorithm \eqref{MUR} for details of this algorithm). Unfortunately, the algorithm may fail to arrive at a stationary point to the problem  \cite{gonza2005,lin2007_1}. A modified multiplicative update algorithm that resolves this issue has been reported \cite[Algorithm 2]{lin2007} along with a theorem \cite[Theorem 7] {lin2007} stating that any accumulation point of the sequence generated by the modified algorithm belongs to the set of all stationary points to Problem \eqref{prob:0}. 
The fast {\tt 'mult'} algorithm, which is based on the multiplicative update algorithms, was implemented in the Statistics and Machine Learning Toolbox in MATLAB R2016a and is given as nnmf functions in the MathWorks documentation \cite{nmf}. Since the {\tt 'mult'} algorithm is more sensitive to initial points, it is useful for finding an optimal solution to Problem \eqref{prob:0} using many random initial points \cite[Description]{nmf}.

Another algorithm used for solving Problem \eqref{prob:0} is the {\em alternating least-squares algorithm} \cite{paat1994} (see Algorithm \eqref{ALS} for details of this algorithm). While there have been reports of convergence in certain special cases (see \cite[p. 160]{berry2007} and references therein), in general, the alternating least-squares algorithm is not guaranteed to converge to a stationary point to Problem \eqref{prob:0}. An alternating least-squares algorithm using projected gradients \cite[Subsection 4.1]{lin2007_1} and one based on alternating nonnegativity constrained least squares  \cite[Section 3]{kim2008} have been reported. 
The fast {\tt 'als'} algorithm, which is based on the alternating least-squares algorithms, was implemented in the Statistics and Machine Learning Toolbox in MATLAB R2016a and is given as nnmf functions in the MathWorks documentation \cite{nmf}. Although the {\tt 'als'} algorithm converges very quickly, it may converge to a point ranked lower than $R$, indicating that the result may not be optimal \cite[Description]{nmf}. 

A stationary point to Problem \eqref{prob:0} is defined as the solution $X^\star \in \mathbb{R}_+^{M \times R} \times \mathbb{R}_+^{R \times N}$ to the following {\em nonmonotone variational inequality} \cite[Definition 1.1.1]{facc1} for the gradient of $f$ over $\mathbb{R}_+^{M \times R} \times \mathbb{R}_+^{R \times N}$:
\begin{align}\label{vi}
\left\langle X - X^\star, \nabla f \left(X^\star \right) \right\rangle \geq 0
\text{ for all } X \in \mathbb{R}_+^{M \times R} \times \mathbb{R}_+^{R \times N},
\end{align}
where $\langle \cdot, \cdot \rangle$ stands for the inner product of $\mathbb{R}^{M \times R} \times \mathbb{R}^{R \times N}$ (see Problem \ref{prob:1} for details of variational inequality). Any local minimizer of Problem \eqref{prob:0} is guaranteed to belong to the set of solutions to Variational Inequality \eqref{vi} \cite[Subchapter 1.3.1]{facc1}. While many iterative algorithms have been proposed for solving variational inequalities \cite[Chapter 12]{facc2}, this paper focuses on using projection methods for solving variational inequalities.

In this paper, we first show that, given $W \in \mathbb{R}^{M \times R}$ and $\mu > 0$, a mapping $T_{W,\mu} \colon \mathbb{R}^{R \times N} \to \mathbb{R}^{R \times N}$ defined for all $H \in \mathbb{R}^{R \times N}$ by
\begin{align}\label{t}
T_{W,\mu} \left(H \right) := P_{\mathbb{R}_+^{R \times N}}
\left[ H - \mu \nabla_H f \left(W,H \right) \right] 
\end{align} 
satisfies the {\em nonexpansivity} condition (i.e., $\|T_{W,\mu}(H_1) - T_{W,\mu}(H_2)\|_{\mathrm{F}} \leq \|H_1 - H_2\|_{\mathrm{F}}$ for all $H_1, H_2 \in \mathbb{R}^{R \times N}$), where $P_{\mathbb{R}_+^{R \times N}}$ stands for the projection onto $\mathbb{R}_+^{R \times N}$ and $\nabla_H f(W, H)$ is the gradient of $f(W, \cdot)$ at $H \in \mathbb{R}^{R \times N}$. Moreover, the {\em fixed point set} of $T_{W,\mu}$, denoted by $\mathrm{Fix}(T_{W,\mu}):= \{H\in \mathbb{R}^{R \times N} \colon T_{W,\mu}(H) = H \}$, coincides with the set of all minimizers of $f(W,\cdot)$ over $\mathbb{R}_+^{R \times N}$ (Proposition \ref{prop:1}).
The same discussion as the one above shows that, given $H \in \mathbb{R}^{R \times N}$ and $\lambda > 0$, a mapping $S_{H,\lambda} \colon \mathbb{R}^{M \times R} \to \mathbb{R}^{M \times R}$ defined for all $W \in \mathbb{R}^{M \times R}$ by
\begin{align}\label{s}
S_{H,\lambda} \left(W \right) := P_{\mathbb{R}_+^{M \times R}}
\left[ W - \lambda \nabla_W f\left(W, H \right) \right]
\end{align}
is nonexpansive and that $\mathrm{Fix}(S_{H,\lambda})$ coincides with the set of all minimizers of $f(\cdot, H)$ over $\mathbb{R}_+^{M \times R}$ (Proposition \ref{prop:2}). From the fact that a search for fixed points of nonexpansive mappings $T_{W,\mu}$ and $S_{H,\lambda}$ defined as in \eqref{t} and \eqref{s} is equivalent to minimization of convex functions $f(W,\cdot)$ and $f(\cdot,H)$ with nonnegativity constraints, fixed point algorithms using two nonexpansive mappings $T_{W,\mu}$ and $S_{H,\lambda}$ should enable a stationary point to Problem \eqref{prob:0} to be found, i.e., a solution to Variational Inequality \eqref{vi}.

There are several useful fixed point algorithms for nonexpansive mappings \cite[Chapter 5]{b-c}, \cite[Chapters 3--6]{berinde}, \cite{halpern,kra,mann,nakajo2003,wit}. The simplest fixed point algorithm is the {\em Krasnosel'ski\u\i-Mann fixed point algorithm} \cite[Chapter 5]{b-c}, \cite{kra,mann} defined for all $n\in \mathbb{N}$ by 
\begin{align}\label{km} 
x_{n+1} := \alpha_n x_n + \left( 1 - \alpha_n \right) T \left(x_n \right),
\end{align}
where $T \colon \mathbb{R}^m \to \mathbb{R}^m$ is nonexpansive with $\mathrm{Fix}(T) \neq \emptyset$, $x_0 \in \mathbb{R}^m$ is an initial point, and $(\alpha_n)_{n\in \mathbb{N}} \subset [0,1]$. The Krasnosel'ski\u\i-Mann fixed point algorithm \eqref{km} can be applied to real-world problems such as signal recovery \cite{comb2007} and analysis of dynamic systems \cite{bot2015}. The sequence $(x_n)_{n\in \mathbb{N}}$ generated by Algorithm \eqref{km} converges to a fixed point of $T$ under certain assumptions \cite[Theorem 5.14]{b-c}. Thanks to this useful fixed point algorithm, an iterative algorithm can be devised for solving Variational Inequality \eqref{vi} in NMF.
 The proposed algorithm is defined as follows: given an $n$th approximation $(W_n, H_n) \in \mathbb{R}_+^{M \times R} \times \mathbb{R}_+^{R \times N}$ of a solution to Variational Inequality \eqref{vi}, compute $(W_{n+1}, H_{n+1})\in \mathbb{R}_+^{M \times R} \times \mathbb{R}_+^{R \times N}$:
\begin{align}\label{proposed}
\begin{split}
&H_{n+1} := \alpha_n H_n + \left(1-\alpha_n \right) T_{W_n, \mu_n} \left(H_n \right),\\
&W_{n+1} := \beta_n W_n + \left(1 - \beta_n \right) S_{H_{n+1}, \lambda_n} \left(W_n \right),
\end{split}
\end{align}
where $T_{W_n,\mu_n}$ and $S_{H_{n+1},\lambda_n}$ are defined as in \eqref{t} and \eqref{s},
and $(\alpha_n)_{n\in \mathbb{N}}, (\beta_n)_{n\in \mathbb{N}} \subset [0,1]$ (see Algorithm \ref{algo:2} for more details).

The first step in Algorithm \eqref{proposed}  is a search for a fixed point of $T_{W_n,\mu_n}$, i.e, a minimizer of $f(W_n,\cdot)$ over $\mathbb{R}_+^{R \times N}$. This step is based on the Krasnosel'ski\u\i-Mann fixed point algorithm. 
The second step in Algorithm \eqref{proposed} is a search for a fixed point of $S_{H_{n+1},\lambda_n}$, i.e., a minimizer of $f(\cdot,H_{n+1})$ over $\mathbb{R}_+^{M \times R}$. The definitions of \eqref{t} and \eqref{s} imply that, for all $n\in \mathbb{N}$, $(W_{n}, H_{n})$ in the proposed algorithm is in $\mathbb{R}_+^{M \times R} \times \mathbb{R}_+^{R \times N}$. Hence, it can be seen intuitively that, for a large enough $n$, the proposed algorithm optimizes not only $f(W_n,\cdot)$ but also $f(\cdot,H_n)$ with the nonnegativity constraints.

One contribution of this work is analysis of the proposed algorithm's convergence. As the second and third paragraphs of this section pointed out, the multiplicative update algorithm and the alternating least-squares algorithm do not always converge to a solution to Variational Inequality \eqref{vi}. In contrast, it is shown that, under the nonempty assumption (Assumption \ref{assum:0}) of the fixed point sets of $T_{W_n,\lambda_n}$ and $S_{H_{n+1},\mu_n}$ $(n\in \mathbb{N})$, any accumulation point of the sequence $((W_n,H_n))_{n\in \mathbb{N}}$ generated by the proposed algorithm belongs to the set of solutions to Variational Inequality \eqref{vi} (Theorem \ref{thm:2}). Additionally, the relationships between Assumption \ref{assum:0} and the assumptions needed to guarantee convergence of the previous fixed point and gradient algorithms \cite{aoyama2007,kuh1981,yamada2005} are discussed.

Another contribution of this work is provision of examples showing the fast convergence and effectiveness of the proposed algorithm for Problem \eqref{prob:0}. The proposed algorithm is numerically compared with the {\tt 'mult'} and {\tt 'als'}  algorithms \cite{nmf} for Problem \eqref{prob:0} under certain conditions. The results show that, when all the elements of $V$ are nonzero and $MN$ is small, the efficiency of the proposed algorithm is almost the same as that of the {\tt 'als'} algorithm (subsection \ref{subsec:4.1}) and that, when all the elements of $V$ are nonzero and $MN$ is large, the proposed algorithm optimizes objective function $f$ in Problem \eqref{prob:0} better than the {\tt 'mult'} and {\tt 'als'} algorithms (subsection \ref{subsec:4.1}). Moreover, the results show that, when $V$ is sparse, the proposed algorithm converges faster than the {\tt 'mult'} and {\tt 'als'} algorithms (subsection \ref{subsec:4.2}) and that the {\tt 'als'} algorithm converges to a point ranked lower than $R$, which is not optimal (subsection \ref{subsec:4.3}).

This paper is organized as follows. Section \ref{sec:2} provides the mathematical preliminaries and explicitly states the main problem (Problem \ref{prob:1}), which is the nonmonotone variational inequality problem in NMF. Section \ref{sec:3} presents the proposed algorithm (Algorithm \ref{algo:2}) for solving the main problem and describes its convergence property (Theorem \ref{thm:2}). Section \ref{sec:4} considers NMF problems in certain situations and numerically compares the behaviors of the {\tt 'mult'} and {\tt 'als'} algorithms  with that of the proposed one. Section \ref{sec:5} concludes the paper with a brief summary and mentions future directions for improving the proposed algorithm.

\section{Mathematical preliminaries}\label{sec:2}
\subsection{Notation and definitions}\label{subsec:2.1}
Let $X^\top$ denote the transpose of a matrix $X$, and let $O$ denote the zero matrix. Let $\mathbb{N}$ denote the set of all positive integers including zero. Let $\mathbb{R}^m$ be an $m$-dimensional Euclidean space with the standard Euclidean norm $\|\cdot\|_2$, and let $\mathbb{R}_+^m := \{(x_1,x_2,\ldots,x_m)^\top \in \mathbb{R}^m \colon x_i \geq 0 \text{ } (i=1,2,\ldots,m)\}$. A mapping $\mathrm{vec}(\cdot) \colon \mathbb{R}^{m\times n} \to \mathbb{R}^{mn}$ is defined for all $X := [X_{ij}] \in \mathbb{R}^{m\times n}$ by $\mathrm{vec}(X) := (X_{11},X_{12},\ldots,X_{1n}, \ldots, X_{m1},X_{m2}\ldots,X_{mn})^\top \in \mathbb{R}^{mn}$. The inner product of two matrices $X,Y\in \mathbb{R}^{m\times n}$ is defined by $X \bullet Y := \mathrm{vec}(X)^\top \mathrm{vec}(Y)$ and induces the Frobenius norm $\|X\|_{\mathrm{F}} := \sqrt{X \bullet X}$. The inner product of $\mathbb{R}^{m_1 \times n_1} \times \mathbb{R}^{m_2 \times n_2}$ is defined for all $(X_1, Y_1), (X_2, Y_2) \in \mathbb{R}^{m_1 \times n_1} \times \mathbb{R}^{m_2 \times n_2}$ by $\langle (X_1, Y_1), (X_2, Y_2) \rangle := X_1 \bullet X_2 + Y_1 \bullet Y_2$.
 
The fixed point set of a mapping $T \colon \mathbb{R}^m \to \mathbb{R}^m$ is denoted by $\mathrm{Fix}(T) := \{x\in \mathbb{R}^m \colon T(x)=x \}$. $T \colon \mathbb{R}^m \to \mathbb{R}^m$ is said to be {\em nonexpansive} \cite[Definition 4.1(ii)]{b-c} if $\|T(x) - T(y)\|_2 \leq \|x-y\|_2$ for all $x,y\in \mathbb{R}^m$. $T$ is said to be {\em firmly nonexpansive} \cite[Definition 4.1(i)]{b-c} if $\|T(x) - T(y)\|_2^2 + \| (\mathrm{Id} - T)(x) - (\mathrm{Id}-T)(y) \|_2^2 \leq \|x-y\|_2^2$ for all $x,y\in \mathbb{R}^m$, where $\mathrm{Id}$ stands for the identity mapping. A firm nonexpansivity condition implies a nonexpansivity condition. The metric projection \cite[Subchapter 4.2, Chapter 28]{b-c} onto a nonempty, closed convex set $C \subset \mathbb{R}^m$, denoted by $P_C$, is defined for all $x\in \mathbb{R}^m$ by $P_C(x) \in C$ and $\|x - P_C(x)\|_2 = \inf_{y\in C}\|x-y\|_2$. Let $x,p\in \mathbb{R}^m$. Then $p = P_C (x)$ if and only if $p\in C$ and $(y-p)^\top (x-p) \leq 0$ for all $y\in C$ \cite[Theorem 3.14]{b-c}. $P_C$ is firmly nonexpansive with $\mathrm{Fix}(P_C) = C$ \cite[Proposition 4.8, (4.8)]{b-c}.

Suppose that $f\colon \mathbb{R}^m \to \mathbb{R}$ is convex and that 
$\nabla f \colon \mathbb{R}^m \to \mathbb{R}^m$ is Lipschitz continuous with the Lipschitz constant $1/L > 0$, i.e., 
$\|\nabla f (x) - \nabla f(y) \|_2 \leq (1/L) \|x-y\|_2$ for all $x,y\in \mathbb{R}^m$. 
Then $\nabla f$ is inverse-strongly monotone with a constant $L$ \cite[Th\'eor\`eme 5]{baillon1977}, i.e.,
$(x-y)^\top (\nabla f (x) - \nabla f(y)) \geq L \|\nabla f(x)-\nabla f(y)\|_2^2$ for all $x,y\in \mathbb{R}^m$.
The set of all minimizers of a function $f\colon \mathbb{R}^m \to \mathbb{R}$ over a set $C \subset \mathbb{R}^m$
is denoted by $\argmin_{x\in C} f(x)$.

\subsection{Nonmonotone variational inequality in NMF}\label{subsec:2.2}
Our main objective here is to find a stationary point to the nonconvex optimization problem \eqref{prob:0} that is defined by a solution to the following {\em nonmonotone variational inequality} \cite[Definition 1.1.1]{facc1} for $\nabla f$ over $\mathbb{R}_+^{M \times R} \times \mathbb{R}_+^{R \times N}$:
\begin{prob}\label{prob:1}
Find a point $X^\star := (W^\star, H^\star)$ in $\mathrm{VI} (C, \nabla f)$ where,
\begin{align*} 
\text{for all } X &:= \left(W,H\right) \in \mathbb{R}^{M \times R} \times \mathbb{R}^{R \times N},\\
\nabla f \left( X \right) 
&= \left( \nabla_W f \left(W, H \right), \nabla_H f \left(W, H \right) \right)\\
&= \left( \left(WH - V \right) H^\top, W^\top \left( WH - V\right)  \right),\\
\mathrm{VI} (C, \nabla f)
&:= \left\{ X^\star \in C := \mathbb{R}_+^{M \times R} \times \mathbb{R}_+^{R \times N} \colon 
 \left\langle X - X^\star, \nabla f \left(X^\star \right) \right\rangle \geq 0 \text{ } 
  \left(X \in C \right) \right\}.
\end{align*}
\end{prob}
Problem \ref{prob:1} is called the {\em stationary point problem} \cite[Subchapter 1.3.1]{facc1} associated with Problem \eqref{prob:0}.
The closedness and convexity conditions of $C:= \mathbb{R}_+^{M \times R} \times \mathbb{R}_+^{R \times N}$ 
guarantee that any local minimizer of Problem \eqref{prob:0} belongs to $\mathrm{VI}(C,\nabla f)$ in Problem \ref{prob:1}. 

First, the following proposition is proven.

\begin{prop}\label{prop:1}
Given $W \in \mathbb{R}^{M \times R}$, let $T_{W,\mu} \colon \mathbb{R}^{R \times N} \to \mathbb{R}^{R \times N}$
be a mapping defined for all $H \in \mathbb{R}^{R \times N}$ by
\begin{align}\label{S}
\begin{split}
T_{W,\mu} \left(H \right) 
&:= P_{\mathbb{R}_+^{R \times N}} 
\left[ H - \mu \nabla_H f \left( W, H\right) \right]\\
&= P_{\mathbb{R}_+^{R \times N}} 
\left[ H - \mu W^\top \left( WH - V\right) \right],
\end{split}
\end{align}
where $\mu \in [0,\infty)$ is defined by
\begin{align}\label{mu}
\mu \in
\begin{cases}
\displaystyle{\left[0, \frac{2}{\left\|W^\top W \right\|_{\mathrm{F}}} \right]} \text{ } &\left(W \neq O \right),\\
\displaystyle{\left[0, \infty \right)} \text{ } &\left(W = O \right).
\end{cases}
\end{align}
Then $T_{W,\mu}$ satisfies the nonexpansivity condition, and
$\mathrm{Fix}(T_{W,\mu}) = \argmin_{H \in \mathbb{R}_+^{R \times N}} f(W,H)$ holds.
\end{prop}

{\em Proof:}
Fix $W$ $(\neq O) \in \mathbb{R}^{M \times R}$ arbitrarily, and let $H_1, H_2 \in \mathbb{R}^{R \times N}$.
The nonexpansivity condition of $P_{\mathbb{R}_+^{R \times N}}$
ensures that 
\begin{align*}
\left\| T_{W,\mu} \left(H_1 \right) - T_{W,\mu} \left(H_2 \right) \right\|_{\mathrm{F}}^2
\leq
\left\| \left( H_1 - \mu \nabla_H f \left( W, H_1\right)  \right) 
      - \left( H_2 - \mu \nabla_H f \left( W, H_2\right)  \right) \right\|_{\mathrm{F}}^2,
\end{align*}
which, together with $\| X - Y \|_{\mathrm{F}}^2 = \| X \|_{\mathrm{F}}^2 - 2 X \bullet Y + \| Y \|_{\mathrm{F}}^2$ $(X,Y \in \mathbb{R}^{m\times n})$, implies that 
\begin{align}\label{ineq:1}
\begin{split}
&\quad \left\| T_{W,\mu} \left(H_1 \right) - T_{W,\mu} \left(H_2 \right) \right\|_{\mathrm{F}}^2\\
&\leq
\left\| H_1 - H_2 \right\|_{\mathrm{F}}^2
-2 \mu \left( H_1 - H_2 \right) \bullet \left(\nabla_H f \left( W, H_1 \right)  - \nabla_H f \left( W, H_2 \right) \right)\\ 
&\quad + \mu^2 \left\|\nabla_H f \left( W, H_1 \right)  - \nabla_H f \left( W, H_2 \right) \right\|_{\mathrm{F}}^2.
\end{split}
\end{align}
The definition of $\nabla_H f ( W, \cdot)$ means that 
\begin{align*}
\left\|\nabla_H f \left( W, H_1 \right) - \nabla_H f \left( W,H_2 \right) \right\|_{\mathrm{F}}
&= 
\left\| W^\top \left(W H_1 - V \right)  - W^\top \left(W H_2 - V \right)  \right\|_{\mathrm{F}}\\
&\leq
\left\|W^\top W \right\|_{\mathrm{F}} \left\| H_1 - H_2 \right\|_{\mathrm{F}},
\end{align*}
which implies that $\nabla_H f (W, \cdot)$ is Lipschitz continuous with 
the Lipschitz constant $\|W^\top W \|_{\mathrm{F}}$.
Accordingly, $\nabla_H f ( W, \cdot)$ is inverse-strongly monotone with the constant $1/\|W^\top W\|_{\mathrm{F}}$ \cite[Th\'eor\`eme 5]{baillon1977} (see subsection \ref{subsec:2.1}).
Hence, \eqref{ineq:1} implies that
\begin{align*}
&\quad \left\| T_{W,\mu} \left(H_1 \right) - T_{W,\mu} \left(H_2 \right) \right\|_{\mathrm{F}}^2\\
&\leq
\left\| H_1 - H_2 \right\|_{\mathrm{F}}^2
+ \mu \left( \mu - \frac{2}{\left\|W^\top W \right\|_{\mathrm{F}}} \right)
\left\|\nabla_H f \left(W, H_1 \right)  - \nabla_H f \left(W, H_2 \right)  \right\|_{\mathrm{F}}^2,
\end{align*}
which, together with \eqref{mu}, means that 
\begin{align*}
\left\| T_{W,\mu} \left(H_1 \right) - T_{W,\mu} \left(H_2 \right) \right\|_{\mathrm{F}}
\leq
\left\| H_1 - H_2 \right\|_{\mathrm{F}}.
\end{align*}
That is, $T_{W,\mu}$ is nonexpansive.
Consider the case where $W = O$. 
Since \eqref{S} and \eqref{mu} ensure that, for all $H\in \mathbb{R}^{R \times N}$ and for all $\mu > 0$, 
\begin{align}\label{zero}
T_{W,\mu} \left(H \right) 
= 
P_{\mathbb{R}_+^{R \times N}} \left[ H - \mu O^\top \left( OH -V \right)\right]
= 
P_{\mathbb{R}_+^{R \times N}} \left(H \right),
\end{align}
$T_{W,\mu}$ is nonexpansive. 

Let $W \in \mathbb{R}^{M \times R}$ and $\mu > 0$. 
Theorem 3.14 in \cite{b-c} (see subsection \ref{subsec:2.1}) implies that $H = T_{W,\mu}(H)$ if and only if 
$(\bar{H} - H) \bullet \nabla_H f(W,H) \geq 0$ for all $\bar{H} \in \mathbb{R}_+^{R \times N}$.
The convexity of $f(W,\cdot)$ thus leads to $\mathrm{Fix}(T_{W,\mu}) = \argmin_{H \in \mathbb{R}_+^{R \times N}}f(W,H)$ \cite[Subchapter 1.3.1]{facc1}.
This completes the proof.
$\Box$

A discussion similar to the proof of Proposition \ref{prop:1}
leads to the following.

\begin{prop}\label{prop:2}
Given $H \in \mathbb{R}^{R \times N}$, let $S_{H,\lambda} \colon \mathbb{R}^{M \times R} \to \mathbb{R}^{M \times R}$
be a mapping defined for all $W \in \mathbb{R}^{M \times R}$ by
\begin{align}\label{T}
\begin{split}
S_{H,\lambda} \left(W \right) 
&:= P_{\mathbb{R}_+^{M \times R}} 
\left[ W - \lambda \nabla_W f \left( W, H \right)  \right]\\
&= P_{\mathbb{R}_+^{M \times R}} 
\left[ W - \lambda \left(WH - V \right) H^\top \right],
\end{split}
\end{align}
where $\lambda \in [0,\infty)$ is defined by
\begin{align}\label{lambda}
\lambda \in
\begin{cases}
\displaystyle{\left[0, \frac{2}{\left\|H^\top H \right\|_{\mathrm{F}}} \right]} \text{ } &\left(H \neq O \right),\\
\displaystyle{\left[0, \infty \right)} \text{ } &\left(H = O \right).
\end{cases}
\end{align}
Then $S_{H,\lambda}$ satisfies the nonexpansivity condition, and 
$\mathrm{Fix}(S_{H,\lambda}) = \argmin_{W \in \mathbb{R}_+^{M \times R}} f(W,H)$ holds.
\end{prop}

The following describes the relationship between $\mathrm{VI}(C,\nabla f)$ and the fixed point sets of $T_{W,\mu}$ and $S_{H,\lambda}$ in Propositions \ref{prop:1} and \ref{prop:2}.

\begin{prop}\label{prop:3}
Let $\mathrm{VI}(C,\nabla f)$ be the solution set of Problem \ref{prob:1} and 
let $T_{W,\mu}$ and $S_{H,\lambda}$ be mappings defined by \eqref{S} and \eqref{T}.
Then
\begin{align*}
\left\{ \left(W,H \right) \in C \colon 
W \in \bigcap_{\substack{H \in \mathbb{R}_+^{R \times N} \\ \lambda > 0}} \mathrm{Fix}\left(S_{H,\lambda}\right), \text{ }
H \in \bigcap_{\substack{W \in \mathbb{R}_+^{M \times R} \\ \mu > 0}} \mathrm{Fix}\left(T_{W,\mu}\right)
 \right\} 
\subset \mathrm{VI}\left(C,\nabla f \right).
\end{align*}
\end{prop}

{\em Proof:}
Let $(\bar{W},\bar{H}) \in D := \{(W,H ) \in C \colon W \in \bigcap_{H \in \mathbb{R}_+^{R \times N},\lambda > 0} \mathrm{Fix}(S_{H,\lambda}), H \in \bigcap_{W \in \mathbb{R}_+^{M \times R},\mu > 0} \mathrm{Fix}(T_{W,\mu}) \}$. Then, for all $\lambda, \mu > 0$, $\bar{W} \in \mathrm{Fix}(S_{\bar{H},\lambda})$, and $\bar{H} \in \mathrm{Fix}(T_{\bar{W},\mu})$. Hence, for all $\lambda, \mu > 0$, $\bar{W} = P_{\mathbb{R}_+^{M \times R}} [\bar{W} - \lambda \nabla_W f(\bar{W}, \bar{H})]$, and $\bar{H} = P_{\mathbb{R}_+^{R \times N}} [\bar{H} - \mu \nabla_H f(\bar{W}, \bar{H})]$, which, together with Theorem 3.14 in \cite{b-c} (see also subsection \ref{subsec:2.1}), guarantees that, for all $\lambda, \mu > 0$, for all $W \in \mathbb{R}_+^{M \times R}$, and for all $H \in \mathbb{R}_+^{R \times N}$, $(W - \bar{W}) \bullet \lambda \nabla_W f(\bar{W}, \bar{H}) \geq 0$ and $(H - \bar{H}) \bullet \mu \nabla_H f(\bar{W}, \bar{H}) \geq 0$, i.e., $D \subset \mathrm{VI}(C,\nabla f)$.
$\Box$

\section{Fixed point algorithm}\label{sec:3}
Several useful algorithms have been presented for solving fixed point problems for nonexpansive mappings, such as the Krasnosel'ski\u\i-Mann fixed point algorithm \cite{kra,mann}, the Halpern fixed point algorithm \cite{halpern,wit}, and the hybrid method \cite{nakajo2003}. This section describes our proposed algorithm based on the Krasnosel'ski\u\i-Mann fixed point algorithm \eqref{km} for solving Problem \ref{prob:1}. For convenience, we rewrite the Krasnosel'ski\u\i-Mann fixed point algorithm as follows: given $x_0 \in \mathbb{R}^m$ and a nonexpansive mapping $T \colon \mathbb{R}^m \to \mathbb{R}^m$,
\begin{align}\label{KM}
x_{n+1} := \alpha_n x_n + \left( 1 - \alpha_n \right) T(x_n) \quad (n\in \mathbb{N}),
\end{align}
where $(\alpha_n)_{n\in \mathbb{N}} \subset [0,1]$ and 
the nonempty condition of $\mathrm{Fix}(T)$ is assumed.
Algorithm \eqref{KM} converges to a fixed point of $T$ if $(\alpha_n)_{n\in \mathbb{N}}$ satisfies 
$\sum_{n=0}^\infty \alpha_n (1 - \alpha_n) = \infty$ \cite[Theorem 5.14]{b-c}.

Algorithm \ref{algo:2} is the proposed algorithm based on Algorithm \eqref{KM} for solving Problem \ref{prob:1}.

\begin{algorithm}                      
\caption{Fixed point algorithm for NMF}         
\label{algo:2}                          
\begin{algorithmic}[1]                  
\REQUIRE $n \in \mathbb{N}$, $V \in \mathbb{R}_+^{M \times N}$, $1 \leq R <  \min\{M,N\}$, $(\alpha_n)_{n\in \mathbb{N}}$, $(\beta_n)_{n\in \mathbb{N}} \subset [0,1]$
\STATE $n \gets 0$, $H_0 \in \mathbb{R}_+^{R \times N}$, $W_0 \in \mathbb{R}_+^{M \times R}$
\REPEAT  
 \IF{$W_n \neq O$} 
  \STATE 
   $\mu_n \in 
   \displaystyle{\left(0, 
   \frac{2}{\left\|W_n^\top W_n\right\|_{\mathrm{F}}} \right]}$
  \STATE
   $T_n \left(H_n \right)
   := P_{\mathbb{R}_+^{R \times N}} 
\left[ H_n - \mu_n W_n^\top \left( W_n H_n - V \right) \right]$ 
 \ELSE
  \STATE 
  $T_n \left(H_n \right) 
  :=  H_n$
  \ENDIF
\STATE $H_{n+1} := \alpha_n H_n + \left( 1 - \alpha_n \right) T_n \left( H_n \right)$\\
 \IF{$H_{n+1} \neq O$} 
  \STATE 
   $\lambda_n \in 
   \displaystyle{\left(0, 
   \frac{2}{\left\|H_{n+1}^\top H_{n+1} \right\|_{\mathrm{F}}} \right]}$
  \STATE
   $S_n \left(W_n \right)
   := P_{\mathbb{R}_+^{M \times R}} 
   \left[ W_n - \lambda_n \left(W_n H_{n+1} - V \right)H_{n+1}^\top \right]$ 
 \ELSE
  \STATE 
  $S_n \left(W_n \right) 
  := W_n$
  \ENDIF
 \STATE 
 $W_{n+1} := \beta_n W_n + \left( 1 - \beta_n \right) S_n \left( W_n \right)$
 \STATE 
 $n \gets n+1$
\UNTIL{stopping condition is satisfied}
\end{algorithmic}
\end{algorithm}

Steps 3--8 in Algorithm \ref{algo:2} generate a nonexpansive mapping $T_n$ on the basis of the results of Proposition \ref{prop:1}. 
When $W_n = O$, \eqref{S} implies that, for all $H \in \mathbb{R}_+^{R \times N}$, $T_{W_n,\mu_n}(H) = P_{\mathbb{R}_+^{R \times N}}(H) = H$ (see also \eqref{zero}), which leads to step 7. Step 9 in Algorithm \ref{algo:2} generates the $(n+1)$th iteration $H_{n+1}$ using the convex combination of $H_n$ and $T_n (H_n)$. This shows that step 9 is based on Algorithm \eqref{KM}. Steps 10--15 generate a nonexpansive mapping $S_n$ on the basis of the results of Proposition \ref{prop:2}. 
This mapping 
depends on $H_{n+1}$ given in step 9. Iteration $W_{n+1}$ is generated by Algorithm \eqref{KM}. The stopping condition in step 18 is, for example, that the number of iterations has reached a positive integer $n_{\max}$ or that $f(W_n, H_n) - f(W_{n-1}, H_{n-1}) < \epsilon$, where $\epsilon > 0$ is small enough. See section \ref{sec:4} for the stopping conditions used in the experiments for the {\tt 'mult'} and {\tt'als'} algorithms \cite{nmf} and Algorithm \ref{algo:2}. 

Analysis of previous fixed point algorithms \cite[Chapter 5]{b-c}, \cite[Chapters 3--5]{berinde}, \cite{aoyama2007_1,aoyama2007,halpern,iiduka_siopt2013,iiduka_mp2015,iiduka_hishinuma_siopt2014,kra,kuh1981,mann,nakajo2003,wit} was based on the nonempty assumption of fixed point sets. Hence, the following was assumed in our analysis of Algorithm \ref{algo:2}.

\begin{assum}\label{assum:0}
Let $S_n$ and $T_n$ ($n\in \mathbb{N}$) be nonexpansive mappings in Algorithm \ref{algo:2}.
Then there exists $m_0 \in \mathbb{N}$ such that 
\begin{align*}
\bigcap_{n=m_0}^\infty \mathrm{Fix}\left(S_n \right) \neq \emptyset \text{ and }
\bigcap_{n=m_0}^\infty \mathrm{Fix}\left(T_n \right) \neq \emptyset.
\end{align*}
\end{assum}

Several iterative methods \cite{aoyama2007,kuh1981} based on Algorithm \eqref{KM} have been proposed for finding a common fixed point of a family of nonexpansive mappings defined on a Banach space under the nonempty condition for the intersection of their mappings (see \cite{aoyama2007_1} for a fixed point algorithm based on the Halpern fixed point algorithm \cite{halpern,wit}). The following method was proposed \cite{kuh1981} for finding a common fixed point of a finite number of nonexpansive mappings. Given a family of nonexpansive mappings $(T_i)_{i=1}^k$ and $\alpha \in (0,1)$,
\begin{align}\label{kuh}
x_{n+1} := \left(1-\alpha \right) x_n + \alpha T_k U_{k-1} \left(x_n \right) \text{ } \left(n\in \mathbb{N} \right),
\end{align}
where $x_0$ is an initial point and $(U_i)_{i=0}^k$ is defined by $U_0 := \mathrm{Id}$, $U_1 := (1-\alpha)\mathrm{Id} + \alpha T_1 U_0, \ldots, U_k := (1-\alpha)\mathrm{Id} + \alpha T_k U_{k-1}$. Theorem 1 in \cite{kuh1981} indicates the strong convergence of the sequence $(x_n)_{n\in \mathbb{N}}$ generated by 
Algorithm \eqref{kuh} to a point in $\bigcap_{i=1}^k \mathrm{Fix}(T_i)$ under the nonempty assumption of $\bigcap_{i=1}^k \mathrm{Fix}(T_i)$. The sequence $(x_n)_{n\in \mathbb{N}}$ is generated using the method in \cite{aoyama2007} for finding a common fixed point of a family of nonexpansive mappings $(T_i)_{i=1}^\infty$ as follows:
\begin{align}\label{ao}
x_{n+1} := \alpha_n x_n + \left(1-\alpha_n \right) \sum_{k=0}^n \beta_k \prod_{i=k+1}^n \left(1-\beta_i \right) T_k \left(x_n \right),
\end{align}
where $(\alpha_n)_{n\in \mathbb{N}}\subset [a,b] \subset (0,1)$ for some $a,b>0$ and $(\beta_n)_{n\in \mathbb{N}} \subset (0,1)$ with $\sum_{n=0}^\infty \beta_n < \infty$. Theorem 3.3 in \cite{aoyama2007} shows the weak convergence of the sequence $(x_n)_{n\in \mathbb{N}}$ generated by Algorithm \eqref{ao} to a point in $\bigcap_{i=1}^\infty \mathrm{Fix}(T_i)$ under the nonempty assumption of $\bigcap_{i=1}^\infty \mathrm{Fix}(T_i)$. Meanwhile, Algorithm \ref{algo:2} generates the sequence $((W_n, H_n))_{n\in \mathbb{N}}$ defined by 
\begin{align*}
&H_{n+1} := \alpha_n H_n + \left( 1 - \alpha_n \right) T_n \left(H_n \right),\\
&W_{n+1} := \beta_n W_n + \left( 1 - \beta_n \right) S_n \left(W_n \right).
\end{align*}
Algorithm \ref{algo:2}, which uses two nonexpansive mappings $T_n$ and $S_n$, is simpler than Algorithms \eqref{kuh} and \eqref{ao}, which use multiple nonexpansive mappings at each iteration.

The following sequence $(u_n)_{n\in \mathbb{N}}$ is generated using the adaptive projected subgradient method \cite[Algorithm (11)]{yamada2005} for solving an asymptotic problem of minimizing a sequence of continuous, convex functions 
$\Theta_n \colon H \to [0,\infty)$ $(n\in \mathbb{N})$ over a nonempty, closed convex set $K \subset H$:
\begin{align}\label{yamada}
u_{n+1} := 
\begin{cases}
\displaystyle{P_K \left[ u_n - \lambda_n \frac{\Theta_n (u_n)}{\| \Theta_n^{'}(u_n) \|^2} \Theta_n^{'}(u_n) \right]}
\text{ } &\left(\Theta_n^{'}(u_n) \neq 0 \right),\\
u_n \text{ } &\left(\Theta_n^{'}(u_n) = 0 \right),
\end{cases}
\end{align}
where $\|\cdot\|$ is the norm of a real Hilbert space $H$, $(\lambda_n)_{n\in \mathbb{N}} \subset [a,b] \subset (0,2)$ for some $a,b > 0$, and $\Theta_n^{'}(u_n)$ stands for the subgradient of $\Theta_n$ at $u_n$. Under the assumptions \cite[Assumptions (12) and (14)]{yamada2005} that $(\Theta_n^{'}(u_n))_{n\in \mathbb{N}}$ is bounded and there exists $n_0 \in \mathbb{N}$ such that
\begin{align}\label{cap0}
\bigcap_{n = n_0}^\infty \argmin_{u\in K} \Theta_n (u) \neq \emptyset \text{ and } 
\inf_{u\in K} \Theta_n (u) = 0 \text{ for all } n \geq n_0,
\end{align}
the sequence $(u_n)_{n\in \mathbb{N}}$ generated by Algorithm \eqref{yamada} satisfies $\lim_{n\to \infty} \Theta_n (u_n) = 0$ \cite[Theorem 2(b)]{yamada2005}. Moreover, if $\bigcap_{n = n_0}^\infty \argmin_{u\in K} \Theta_n (u)$ has an interior point, $(u_n)_{n\in \mathbb{N}}$ in Algorithm \eqref{yamada} strongly converges to $u^\star \in K$ with $\lim_{n\to \infty} \Theta_n(u^\star) = 0$ \cite[Theorem 2(c)]{yamada2005}. Section 4 in \cite{yamada2005} provides some examples of satisfying Assumption \eqref{cap0} and shows that, if there is no interior point of $\bigcap_{n = n_0}^\infty \argmin_{u\in K} \Theta_n (u)$, the normalized least squares algorithm \cite[Chapter 5]{haykin2002}, which is an example of Algorithm \eqref{yamada}, diverges \cite[Example 1(b)]{yamada2005}. Therefore, Assumption \eqref{cap0} and the existence of an interior point of $\bigcap_{n = n_0}^\infty \argmin_{u\in K} \Theta_n (u)$ are crucial to ensuring the convergence of Algorithm \eqref{yamada}.

Similarly to Algorithm \eqref{yamada}, Algorithm \ref{algo:2} is a projected gradient algorithm (see steps 3--8 and 10--15 in Algorithm \ref{algo:2}). Propositions \ref{prop:1} and \ref{prop:2} guarantee that Assumption \ref{assum:0} can be expressed as
\begin{align}\label{assum:0_1}
\bigcap_{n=m_0}^\infty \argmin_{W \in \mathbb{R}_+^{M \times R}} f_n (W) \neq \emptyset
\text{ and }
\bigcap_{n=m_0}^\infty \argmin_{H \in \mathbb{R}_+^{R \times N}} g_n (H) \neq \emptyset,
\end{align}
where, for all $n\in \mathbb{N}$, $f_n$ and $g_n$ are convex functions defined by 
$f_n (W) := f(W, H_{n+1})$ $(W \in \mathbb{R}^{M \times R})$ and $g_n (H) := f(W_n,H)$ $(H \in \mathbb{R}^{R \times N})$.
While both \eqref{yamada} \cite[Algorithm (11)]{yamada2005} and Algorithm \ref{algo:2} work under Assumption \eqref{cap0}, Algorithm \ref{algo:2} can also work under a weaker assumption \eqref{assum:0_1} than \eqref{cap0}.
In the special case where, for all $n \geq m_0$, $\inf_{W \in \mathbb{R}_+^{M \times R}} f_n(W) = 0$ and $\inf_{H \in \mathbb{R}_+^{R \times N}} g_n(H) = 0$ \cite[Assumption (14)]{yamada2005}, Assumption \eqref{assum:0_1} is satisfied if and only if the matrix equations $W_n H = V$ and $W H_{n+1} = V$ are consistent for all $n \geq m_0$, i.e., there exist $\bar{W}$ and $\bar{H}$ such that, for all $n \geq m_0$, $W_n \bar{W} V = V$ and $V \bar{H} H_{n+1} = V$ \cite[Theorem 1, Subchapter 2.1]{ben2003} (see \cite[Theorem 1, Subchapter 2.1]{ben2003} for the closed forms of the general solutions to the matrix equations $W_n H = V$ and $W H_{n+1} = V$ $(n \geq m_0)$).

Previously reported results \cite[Subchapter 4.5]{b-c}, \cite{bell1966,browder1965,bruck1974,demarr1963} provide the properties of fixed point sets of nonexpansive mappings and sufficient conditions for the existence of a common fixed point of a family of nonexpansive mappings.\footnote{These results were obtained under the condition that the whole space is a Banach space or a Hilbert space.} See \cite{iiduka_siopt2013,iiduka_hishinuma_siopt2014,yamada2005} for examples of satisfying the nonempty condition for the intersection of fixed point sets of nonexpansive mappings in signal processing and network resource allocation.

The discussion in this section is based on the following assumption.

\begin{assum}\label{assum:2}
The sequences $(\alpha_n)_{n\in \mathbb{N}}$, $(\beta_n)_{n\in \mathbb{N}}$, $(\mu_n)_{n\in \mathbb{N}}$,
and $(\lambda_n)_{n\in \mathbb{N}}$ in Algorithm \ref{algo:2} satisfy 
\begin{align*}
&\text{{\em (C1)} } 0 < \liminf_{n\to \infty} \alpha_n \leq \limsup_{n\to \infty} \alpha_n < 1,\text{ }
0 < \liminf_{n\to \infty} \beta_n \leq \limsup_{n\to \infty} \beta_n < 1,\\
&\text{{\em (C2)} } 0 < \liminf_{n\to \infty} \mu_n,\text{ }
0 < \liminf_{n\to \infty} \lambda_n.
\end{align*}
\end{assum}

Examples of $(\alpha_n)_{n\in \mathbb{N}}$ and $(\beta_n)_{n\in \mathbb{N}}$ satisfying (C1)
are $\alpha_n := \alpha \in (0,1)$ and $\beta_n := \beta \in (0,1)$ $(n\in \mathbb{N})$.
Under Assumptions \ref{assum:0} and (C1), the boundedness conditions of $(W_n)_{n\in \mathbb{N}}$
and $(H_n)_{n\in \mathbb{N}}$ are guaranteed (see Lemma \ref{lem:3_2}(i)).
Accordingly, $\mu_n$ and $\lambda_n$ $(n\in \mathbb{N})$ can be chosen such that
\begin{align}\label{para}
\mu_n := \frac{2}{\left\| W_n^\top W_n \right\|_{\mathrm{F}}} \text{ and } 
\lambda_n := \frac{2}{\left\| H_{n+1}^\top H_{n+1} \right\|_{\mathrm{F}}},
\end{align}
which satisfy Assumption (C2).

The following is a convergence analysis of Algorithm \ref{algo:2}.

\begin{thm}\label{thm:2}
Under Assumptions \ref{assum:0} and \ref{assum:2}, any accumulation point of the sequence $((W_n,H_n))_{n\in \mathbb{N}}$
generated by Algorithm \ref{algo:2} belongs to $\mathrm{VI}(C,\nabla f)$.
\end{thm}

Since Assumptions \ref{assum:0} and (C1) ensure the boundedness of $((W_n,H_n))_{n\in \mathbb{N}}$ in Algorithm \ref{algo:2} (see Lemma \ref{lem:3_2}), there exists an accumulation point of $((W_n,H_n))_{n\in \mathbb{N}}$. Accordingly, Theorem \ref{thm:2} implies the existence of a point in $\mathrm{VI}(C,\nabla f)$. In the case where $\mathrm{VI}(C,\nabla f)$ consists of one point, Theorem \ref{thm:2} guarantees that the whole sequence $((W_n,H_n))_{n\in \mathbb{N}}$ generated by Algorithm \ref{algo:2} converges to that point.

Let us compare Algorithm \ref{algo:2} with previous algorithms for solving Problem \eqref{prob:0}. A useful approach to solving Problem \eqref{prob:0} is the following multiplicative update algorithm \cite[Chapter 3]{cic2009}, \cite[Algorithm (4)]{lee2001}, \cite[Algorithm 1]{lin2007} (see \cite[Description, {\tt 'mult'} algorithm]{nmf}, \cite[p. 158]{berry2007}, and \cite[Algorithm 2]{lin2007} for modified multiplicative update algorithms): given $H_n := [H_{n,bj}] = [(H_n)_{bj}]$ and $W_n := [W_{n,ia}] = [(W_n)_{ia}]$,
\begin{align}\label{MUR}
\begin{split}
&H_{n+1,bj} := H_{n,bj} \frac{((W_n)^\top V)_{bj}}{(W_n^\top W_n H_n)_{bj}} \text{ } (b = 1,2,\ldots,R, \text{ }
j= 1,2,\ldots,N),\\
&W_{n+1,ia} := W_{n,ia} \frac{(VH_{n+1}^\top)_{ia}}{(W_n H_{n+1} H_{n+1}^\top)_{ia}}
\text{ } (i = 1,2,\ldots,M, \text{ }
a= 1,2,\ldots,R),
\end{split}
\end{align}
which satisfies that $f(W_n, H_{n+1}) \leq f(W_n, H_n)$ and that $f(W_{n+1},H_{n+1}) \leq f(W_n,H_{n+1})$ for all $n\in \mathbb{N}$. However, Section 5 in \cite{gonza2005} and Section 6 in \cite{lin2007_1} show that Algorithm \eqref{MUR} does not always converge to a stationary point to Problem \eqref{prob:0}. The most well-known algorithms for NMF are the alternating least-squares algorithms \cite[Description, {\tt 'als'} algorithm]{nmf}, \cite[Subsection 3.3]{berry2007}, \cite[Chapter 4]{cic2009}, \cite[Subsection 4.1]{lin2007_1}, \cite{paat1994}. The framework of the basic alternating least-squares algorithm \cite[p. 160]{berry2007} is as follows: given $W_n \in \mathbb{R}_+^{M \times R}$,
\begin{align}\label{ALS}
\begin{split}
H_{n+1} :=  P_{\mathbb{R}_+^{R \times N}} \left(\bar{H}_{n} \right), \text{ }
W_{n+1} := P_{\mathbb{R}_+^{M\times R}} \left(\bar{W}_{n} \right),
\end{split}
\end{align}
where $\bar{H}_{n}$ satisfies $W_n^\top W_n \bar{H}_{n} = W_n^\top V$ and $\bar{W}_{n}$ satisfies $H_{n+1} H_{n+1}^\top \bar{W}_{n}^\top = H_{n+1} V^\top$. While alternating least-squares algorithms have been shown to converge very quickly, there has been no reported analysis of their convergence, and there is no guarantee that they converge to stationary points to Problem \eqref{prob:0} \cite[p. 160]{berry2007}.

In contrast, Theorem \ref{thm:2} guarantees the convergence of Algorithm \ref{algo:2} to a stationary point to Problem \eqref{prob:0} under Assumptions \ref{assum:0} and \ref{assum:2}. The numerical comparison in section \ref{sec:4} of Algorithm \ref{algo:2} with the {\tt 'mult'} and {\tt 'als'} algorithms given in MATLAB R2016a shows that Algorithm \ref{algo:2} is effective and converges quickly. It also shows that the {\tt 'als'} algorithm converged to a point ranked lower than $R$, which may indicate that the result is not optimal \cite[Description]{nmf}.

\subsection{Proof of Theorem \ref{thm:2}}\label{subsec:3.1}
The proof of Theorem \ref{thm:2} is divided into three steps. First, we prove the following lemma.

\begin{lem}\label{lem:3_2}
Suppose that Assumptions \ref{assum:0} and (C1) hold.
Then 
\begin{enumerate}
\item[{\em (i)}]
there exist $\lim_{n\to \infty} \| H_n - H\|_{\mathrm{F}}$ for all $H \in \bigcap_{n=m_0}^\infty \mathrm{Fix}(T_n)$ and 
$\lim_{n\to \infty} \| W_n - W\|_{\mathrm{F}}$ for all $W \in \bigcap_{n=m_0}^\infty \mathrm{Fix}(S_n)$;
\item[{\em (ii)}]
$\lim_{n\to\infty} \|H_n - T_n(H_n)\|_{\mathrm{F}} = 0$
and
$\lim_{n\to\infty} \|W_n - S_n (W_n)\|_{\mathrm{F}} = 0$.
\end{enumerate}
\end{lem}

{\em Proof:}
(i)
Fix $H \in \bigcap_{n=m_0}^\infty \mathrm{Fix}(T_n)$ arbitrarily.
From $\| \alpha X + (1-\alpha) Y\|_{\mathrm{F}}^2 
= \alpha \| X\|_{\mathrm{F}}^2 + (1-\alpha) \|Y\|_{\mathrm{F}}^2 - \alpha (1-\alpha)\|X - Y\|_{\mathrm{F}}^2$ $(X,Y \in \mathbb{R}^{m \times n}, \alpha \in \mathbb{R})$,
for all $n\geq m_0$,
\begin{align*}
\left\| H_{n+1} - H \right\|_{\mathrm{F}}^2
&=
\left\| \alpha_n \left(H_{n} - H \right) + \left(1-\alpha_n \right) \left(T_n\left(H_n\right) - H \right) \right\|_{\mathrm{F}}^2\\
&= 
\alpha_n \left\|H_{n} - H \right\|_{\mathrm{F}}^2
+ \left(1-\alpha_n \right) \left\| T_n\left(H_n\right) - T_n\left(H\right)  \right\|_{\mathrm{F}}^2\\
&\quad - \alpha_n \left(1-\alpha_n \right) \left\| H_n - T_n\left(H_n\right)  \right\|_{\mathrm{F}}^2,
\end{align*}
which, together with the nonexpansivity of $T_n$ $(n\in \mathbb{N})$ (see Proposition \ref{prop:1}), implies that, for all $n\geq m_0$,
\begin{align}\label{mono}
\begin{split}
\left\| H_{n+1} - H \right\|_{\mathrm{F}}^2
&\leq 
\left\|H_{n} - H \right\|_{\mathrm{F}}^2
- \alpha_n \left(1-\alpha_n \right) \left\| H_n - T_n\left(H_n\right)  \right\|_{\mathrm{F}}^2\\
&\leq \left\|H_{n} - H \right\|_{\mathrm{F}}^2.
\end{split}
\end{align}
Since $(\|H_n - H\|_{\mathrm{F}})_{n\in \mathbb{N}}$ is bounded below and monotone decreasing, 
there exists $\lim_{n\to\infty} \| H_{n+1} - H\|_{\mathrm{F}}$.
Accordingly, $(H_n)_{n\in \mathbb{N}}$ is bounded.

Fix $W \in \bigcap_{n=m_0}^\infty \mathrm{Fix}(S_n)$ arbitrarily. 
Proposition \ref{prop:2} guarantees that, for all $n\in \mathbb{N}$, $S_n$ is nonexpansive.
Therefore, a discussion similar to the one for obtaining \eqref{mono} ensures that, for all $n\geq m_0$, 
\begin{align}\label{mono_1}
\begin{split}
\left\| W_{n+1} - W \right\|_{\mathrm{F}}^2
&\leq 
\left\|W_{n} - W \right\|_{\mathrm{F}}^2
- \beta_n \left(1-\beta_n \right) \left\| W_n - S_n\left(W_n\right)  \right\|_{\mathrm{F}}^2\\
&\leq \left\|W_{n} - W \right\|_{\mathrm{F}}^2,
\end{split}
\end{align}
which implies the existence of $\lim_{n\to\infty}\| W_n - W\|_{\mathrm{F}}$ and the boundedness of $(W_n)_{n\in \mathbb{N}}$.

(ii)
Inequality \eqref{mono} leads to the finding that 
$\lim_{n\to\infty} \alpha_n (1-\alpha_n ) \| H_n - T_n\left(H_n\right)\|_{\mathrm{F}}^2 = 0$,
which, together with Assumption (C1), implies that 
$\lim_{n\to\infty} \| H_n - T_n (H_n) \|_{\mathrm{F}} = 0$.
Moreover, Assumption (C1) and \eqref{mono_1} lead to $\lim_{n\to\infty} \| W_n - S_n(W_n)\|_{\mathrm{F}} = 0$,
which completes the proof.
$\Box$

Lemma \ref{lem:3_2} leads to the following.

\begin{lem}\label{lem:3_1}
Suppose that Assumptions \ref{assum:0} and \ref{assum:2} hold.
Then for all $(W,H) \in C := \mathbb{R}_+^{M \times R} \times \mathbb{R}_+^{R \times N}$,
\begin{align*}
&\limsup_{n\to\infty}\left(T_n\left(H_n \right)- H \right) \bullet \nabla_H f \left(W_n, H_n \right) \leq 0,\\
&\limsup_{n\to\infty}\left(S_n\left(W_n \right)- W\right) \bullet \nabla_W f \left(W_n, H_{n+1} \right)
\leq 0.
\end{align*}
\end{lem}

{\em Proof:}
Fix $(W,H) \in C$ arbitrarily.
The firm nonexpansivity of $\mathbb{R}_+^{R \times N}$ 
and the condition $H = P_{\mathbb{R}_+^{R \times N}}(H)$ $(H \in \mathbb{R}_+^{R \times N})$ ensure that,
for all $n\in \mathbb{N}$, 
\begin{align*}
\left\| T_n\left(H_n\right) - H \right\|_{\mathrm{F}}^2
&=
\left\| P_{\mathbb{R}_+^{R \times N}} 
\left[ H_n - \mu_n \nabla_H f \left(W_n, H_n\right) 
\right]
-
P_{\mathbb{R}_+^{R \times N}} \left(H\right) \right\|_{\mathrm{F}}^2\\
&\leq
\left\| \left(H_n - \mu_n \nabla_H f \left(W_n, H_n \right) \right) - H \right\|_{\mathrm{F}}^2\\
&\quad -
\left\| \left(H_n - \mu_n \nabla_H f \left(W_n, H_n \right) \right) - T_n\left(H_n\right) \right\|_{\mathrm{F}}^2,
\end{align*}
which, together with $\| X - Y \|_{\mathrm{F}}^2 = \| X \|_{\mathrm{F}}^2 - 2 X \bullet Y + \| Y \|_{\mathrm{F}}^2$ $(X,Y \in \mathbb{R}^{m\times n})$, implies that, for all $n\in \mathbb{N}$, 
\begin{align*}
&\quad \left\| T_n\left(H_n\right) - H \right\|_{\mathrm{F}}^2\\
&\leq
\left\| H_n - H \right\|_{\mathrm{F}}^2
-2 \mu_n \left(H_n - H \right) \bullet  \nabla_H f \left(W_n, H_n \right) \\
&\quad 
- \left\| H_n - T_n\left(H_n\right) \right\|_{\mathrm{F}}^2
+ 2 \mu_n \left(H_n - T_n\left(H_n\right) \right) \bullet  \nabla_H f \left(W_n, H_n \right) \\
&=
\left\| H_n - H \right\|_{\mathrm{F}}^2
- \left\| H_n - T_n\left(H_n\right) \right\|_{\mathrm{F}}^2
+ 2 \mu_n \left(H - T_n\left(H_n\right) \right) \bullet  \nabla_H f \left(W_n, H_n\right).
\end{align*}
Accordingly, the convexity of $\|\cdot\|_{\mathrm{F}}^2$ implies that, for all $n\in \mathbb{N}$, 
\begin{align}\label{liminf_0}
\begin{split}
\left\| H_{n+1} - H \right\|_{\mathrm{F}}^2
&\leq
\alpha_n \left\| H_n - H \right\|_{\mathrm{F}}^2
+ \left( 1 - \alpha_n \right) \left\| T_n\left(H_n\right) - H \right\|_{\mathrm{F}}^2\\
&\leq
\alpha_n \left\| H_n - H \right\|_{\mathrm{F}}^2
+ \left( 1 - \alpha_n \right)
\big\{
\left\| H_n - H \right\|_{\mathrm{F}}^2
- \left\| H_n - T_n\left(H_n\right) \right\|_{\mathrm{F}}^2\\
&\quad 
+ 2 \mu_n \left(H - T_n\left(H_n\right) \right) \bullet  \nabla_H f \left(W_n, H_n \right) \big\}\\
&\leq
\left\| H_n - H \right\|_{\mathrm{F}}^2
 + 2 \mu_n \left(1 - \alpha_n\right) \left(H - T_n\left(H_n\right) \right) \bullet  \nabla_H f \left(W_n, H_n\right),
\end{split} 
\end{align} 
which implies that, for all $n\in \mathbb{N}$,
\begin{align}\label{limsup_0}
\begin{split}
&\quad 2 \mu_n \left(1 - \alpha_n\right) \left(T_n\left(H_n\right) - H\right) \bullet  \nabla_H f \left(W_n, H_n\right) \\
&\leq
\left\| H_n - H \right\|_{\mathrm{F}}^2 - \left\| H_{n+1} - H \right\|_{\mathrm{F}}^2\\
&=
\left(\left\| H_n - H \right\|_{\mathrm{F}} + \left\| H_{n+1} - H \right\|_{\mathrm{F}} \right)
\left(\left\| H_n - H \right\|_{\mathrm{F}} - \left\| H_{n+1} - H \right\|_{\mathrm{F}} \right)\\
&\leq 
\left(\left\| H_n - H \right\|_{\mathrm{F}} + \left\| H_{n+1} - H \right\|_{\mathrm{F}} \right) \left\| H_{n} - H_{n+1} \right\|_{\mathrm{F}}\\
&\leq
M_1 \left\| H_{n+1} - H_n \right\|_{\mathrm{F}},
\end{split}
\end{align}
where 
the second inequality comes from the triangle inequality,
the third inequality comes from $M_1 := \sup \{ \|H_n - H\|_{\mathrm{F}} + \|H_{n+1} - H\|_{\mathrm{F}} \colon n\in \mathbb{N} \}$, and $M_1 < \infty$ is guaranteed from the boundedness of $(H_n)_{n\in \mathbb{N}}$ (see Lemma \ref{lem:3_2}(i)).
Since the definition of $H_{n+1}$ $(n\in \mathbb{N})$ means that, for all $n\in \mathbb{N}$,
$\| H_{n+1} - H_n \|_{\mathrm{F}} = (1-\alpha_n) \| H_{n} - T_n (H_n) \|_{\mathrm{F}}$,
Lemma \ref{lem:3_2}(ii) and (C1) ensure that
\begin{align}\label{H_n}
\lim_{n\to\infty} \left\| H_{n+1} - H_n \right\|_{\mathrm{F}} = 0.
\end{align}
Hence, \eqref{limsup_0} and \eqref{H_n} lead to the deduction that
\begin{align*}
\limsup_{n\to \infty} \mu_n \left(1 - \alpha_n\right) \left(T_n\left(H_n\right) - H\right) \bullet  \nabla_H f \left(W_n,H_n\right) \leq 0.
\end{align*}
Accordingly, for all $\epsilon > 0$, there exists $n_0 \in \mathbb{N}$ such that, for all $n\geq n_0$,
\begin{align}\label{limsup_1}
\mu_n \left(1 - \alpha_n\right) \left(T_n\left(H_n\right) - H\right) \bullet  \nabla_H f \left(W_n, H_n\right) \leq \epsilon.
\end{align}
Since (C1) and (C2) guarantee that there exists $c_1 > 0$ such that, for all $n \in \mathbb{N}$, 
$0 < 1/(\mu_n(1-\alpha_n)) \leq c_1$,
\eqref{limsup_1} leads to the finding that 
\begin{align*}
\limsup_{n\to \infty} \left(T_n\left(H_n\right) - H\right) \bullet  \nabla_H f \left(W_n, H_n\right) \leq c_1 \epsilon,
\end{align*}
which, together with the arbitrary condition of $\epsilon$, implies that 
\begin{align}\label{limsup_2}
\limsup_{n\to \infty} \left(T_n\left(H_n\right) - H\right) \bullet  \nabla_H f \left(W_n, H_n \right) \leq 0.
\end{align}
 
A discussion similar to the one for obtaining \eqref{limsup_0} leads to the finding that 
there exists $M_2 \in \mathbb{R}$ such that, for all $n\in \mathbb{N}$,
\begin{align*}
\quad 2 \lambda_n \left(1 - \beta_n\right) \left(S_n\left(W_n\right) - W\right) \bullet  \nabla_W f \left(W_n, H_{n+1} \right) 
\leq
M_2 \left\| W_{n+1} - W_n \right\|_{\mathrm{F}}.
\end{align*}
Since the same manner of argument as in the proof of \eqref{H_n} leads to 
$\lim_{n\to \infty} \| W_{n+1} - W_n\|_{\mathrm{F}} = 0$,
it can be found that
\begin{align*}
\limsup_{n\to\infty} \lambda_n \left(1 - \beta_n\right) \left(S_n\left(W_n\right) - W\right) \bullet  \nabla_W f \left(W_n, H_{n+1} \right) \leq 0.
\end{align*}
Therefore, a discussion similar to the one for obtaining \eqref{limsup_2}, together with (C1) and (C2), implies that
\begin{align*}
\limsup_{n\to\infty} \left(S_n\left(W_n\right) - W\right) \bullet  \nabla_W f \left(W_n , H_{n+1} \right) \leq 0,
\end{align*}
which completes the proof.
$\Box$

Lemmas \ref{lem:3_2} and \ref{lem:3_1} lead to the following. 

\begin{lem}
Suppose that Assumptions \ref{assum:0} and \ref{assum:2} hold.
Let $(W_\star, H_\star) \in \mathbb{R}^{M \times R} \times \mathbb{R}^{R \times N}$ be an accumulation point of $((W_n, H_n))_{n\in \mathbb{N}}$ generated by Algorithm \ref{algo:2}.
Then $(W_\star, H_\star) \in \mathrm{VI}(C,\nabla f)$.
\end{lem}

{\em Proof:}
Since Lemma \ref{lem:3_2}(i) ensures the boundedness of $((W_n, H_n))_{n\in \mathbb{N}}$,
there exists an accumulation point of $((W_n, H_n))_{n\in \mathbb{N}}$.
Fix an accumulation point of $((W_n, H_n))_{n\in \mathbb{N}}$, denoted by $(W_\star, H_\star) \in \mathbb{R}^{M \times R} \times \mathbb{R}^{R \times N}$, arbitrarily.
Then there exists a subsequence $((W_{n_i}, H_{n_i}))_{i\in \mathbb{N}}$ of $((W_n, H_n))_{n\in \mathbb{N}}$ such that $((W_{n_i}, H_{n_i}))_{i\in \mathbb{N}}$ converges to $(W_\star, H_\star)$.
Accordingly, \eqref{H_n} ensures that $(H_{n_i + 1})_{i\in \mathbb{N}}$ converges to $H_\star$.
From $((W_n, H_n))_{n\in \mathbb{N}} \subset C := \mathbb{R}_+^{M \times R} \times \mathbb{R}_+^{R \times N}$, 
the closedness of $C$ guarantees that
$(W_\star, H_\star) \in C$. 

Meanwhile, for all $i\in \mathbb{N}$,
\begin{align*}
&\quad \left\| \nabla_H f \left(W_{n_i},H_{n_i} \right)  - \nabla_H f \left(W_\star, H_\star \right) \right\|_{\mathrm{F}}\\
&= 
\left\| W_{n_i}^\top \left( W_{n_i} H_{n_i} - V  \right) - W_{\star}^\top \left( W_{\star} H_{\star} - V  \right) \right\|_{\mathrm{F}}\\
&= 
\left\| \left(W_{\star} - W_{n_i}\right)^\top \left(V - W_\star H_\star \right) + W_{n_i}^\top W_{n_i} \left( H_{n_i} - H_\star \right)
+  W_{n_i}^\top \left( W_{n_i} - W_\star \right)H_\star \right\|_{\mathrm{F}},
\end{align*}
which, together with the triangle inequality and the boundedness of $(W_n)_{n\in \mathbb{N}}$,
implies that  
\begin{align}\label{0}
\lim_{i\to \infty} \left\| \nabla_H f \left(W_{n_i},H_{n_i} \right) - \nabla_H f \left(W_\star,H_\star \right) \right\|_{\mathrm{F}} = 0.
\end{align}
A discussion similar to the one for obtaining \eqref{0}, together with the boundedness of $(H_n)_{n\in \mathbb{N}}$ and $\lim_{i\to \infty} \|H_{n_i + 1} - H_\star \|_{\mathrm{F}} = 0$, leads to 
\begin{align}\label{1}
\lim_{i\to \infty} \left\| \nabla_W f \left(W_{n_i}, H_{n_i + 1} \right) - \nabla_W f \left(W_\star, H_\star \right)  \right\|_{\mathrm{F}} = 0.
\end{align}
Moreover, Lemma \ref{lem:3_2}(ii) guarantees that 
\begin{align}\label{2}
\lim_{i\to \infty} \left\| S_{n_i} \left(W_{n_i} \right) - W_\star \right\|_{\mathrm{F}} = 0 \text{ and }
\lim_{i\to \infty} \left\| T_{n_i} \left(H_{n_i} \right) - H_\star \right\|_{\mathrm{F}} = 0.
\end{align}
Lemma \ref{lem:3_1} leads to the deduction that, for all $(W,H) \in C$,
\begin{align*}
&\quad 
\lim_{i\to\infty}\left(T_{n_i}\left(H_{n_i} \right)- H \right) \bullet \nabla_H f \left(W_{n_i}, H_{n_i} \right)\\
&\leq \limsup_{n\to\infty}\left(T_n\left(H_n \right)- H \right) \bullet \nabla_H f \left(W_n, H_n \right) 
\leq 0,\\ 
&\quad \lim_{i\to\infty}\left(S_{n_i}\left(W_{n_i} \right)- W\right) \bullet \nabla_W f \left(W_{n_i}, H_{n_i +1} \right) \\
&\leq \limsup_{n\to\infty}\left(S_n\left(W_n \right)- W\right) \bullet \nabla_W f \left(W_n , H_{n+1} \right) 
\leq 0.
\end{align*}
Hence, \eqref{0}, \eqref{1}, and \eqref{2} guarantee that, for all $(W,H)\in C$,
\begin{align*}
\left(H_\star - H \right) \bullet \nabla_H f \left(W_\star, H_\star \right)  \leq 0 \text{ and }
\left(W_\star - W\right) \bullet \nabla_W f \left(W_\star, H_\star \right) 
\leq 0,
\end{align*}
which implies that $(W_\star, H_\star) \in \mathrm{VI}(C,\nabla f)$.
This completes the proof.
$\Box$
  
\section{Numerical experiments}\label{sec:4}
This section applies Algorithm \ref{algo:2}, the multiplicative update algorithm \cite[{\tt 'mult'} algorithm]{nmf}, and the alternating least-squares algorithm \cite[{\tt 'als'} algorithm]{nmf} to the following optimization problem \cite[Description]{nmf}: given $V \in \mathbb{R}_+^{M \times N}$ and $1 \leq R < \min \{M,N\}$,
\begin{align}\label{nmf}
\text{minimize } 
F\left(W,H \right) := \frac{\left\| V - WH \right\|_{\mathrm{F}}}{\sqrt{MN}} 
\text{ subject to } 
\left(W,H \right) \in C,
\end{align}
where $C := \mathbb{R}_+^{M \times R} \times \mathbb{R}_+^{R \times N}$. A solution $(W^\star, H^\star)$ to Problem \eqref{nmf} needs to satisfy the following conditions:
\begin{enumerate}
\item[(i)] $H^\star$ is normalized so that the rows of $H^\star$ have unit length;
\item[(ii)] The columns of $W^\star$ are ordered by decreasing length.
\end{enumerate}

The experiments were done using a MacBook Air with a 2.20 GHz Intel(R) Core(TM) i7-5650U CPU processor, 8 GB 1600 MHz DDR3 memory, and Mac OS X El Capitan (Version 10.11.6) operating system. The algorithms used were written in MATLAB R2016a (9.0.0.341360):
\begin{itemize}
\item
\textbf{MULT}: the multiplicative update algorithm \cite[{\tt 'mult'} algorithm]{nmf} 
implemented in Statistics and Machine Learning Toolbox in MATLAB R2016a. 
\item
\textbf{ALS}: the alternating least-squares algorithm \cite[{\tt 'als'} algorithm]{nmf} implemented in Statistics and Machine Learning Toolbox in MATLAB R2016a.
\item
\textbf{Proposed ($C$)}: Algorithm \ref{algo:2} when $\alpha_n = \beta_n = C \in (0,1)$ $(n\in \mathbb{N})$,
$\mu_n := 2/\max \{ 1, \|W_n^\top W_n \|_{\mathrm{F}} \}$, and $\lambda_n := 2/\max \{ 1, \|H_{n+1}^\top  H_{n+1} \|_{\mathrm{F}} \}$ $(n\in \mathbb{N})$.
\item
\textbf{Proposed ($C, c$)}: Algorithm \ref{algo:2} when $\alpha_n = \beta_n = C \in (0,1)$ $(n\in \mathbb{N})$ and 
$\mu_n = \lambda_n = c > 0$ $(n\in \mathbb{N})$.
\end{itemize}

Although the MULT and ALS algorithms converge very quickly, their convergence to a stationary point to Problem \eqref{nmf} is not guaranteed. Hence, they may converge to a point ranked lower than $R$, indicating that the result is not optimal \cite[Description]{nmf}. Meanwhile, the convergence of Algorithm \ref{algo:2} to a stationary point to Problem \eqref{nmf} is guaranteed under the assumptions in Theorem \ref{thm:2}. The parameters $\alpha_n$ and $\beta_n$ $(n\in \mathbb{N})$ in Algorithm \ref{algo:2} were set to satisfy Assumption (C1). The parameters $\mu_n := 2/\max \{1, \|W_n^\top W_n \|_{\mathrm{F}} \}$ and $\lambda_n := 2/\max \{1, \|H_{n+1}^\top H_{n+1} \|_{\mathrm{F}} \}$ $(n\in \mathbb{N})$ were used in Proposed $(C)$ to satisfy Assumption (C2) and steps 4 and 11 in Algorithm \ref{algo:2} (see also \eqref{para}).\footnote{Suppose that $W_n \neq O$. If $\|W_n^\top W_n\|_{\mathrm{F}} \leq 1$, $\mu_n := 2/\max \{1, \|W_n^\top W_n \|_{\mathrm{F}} \} = 2 \in (0, 2/\|W_n^\top W_n\|_{\mathrm{F}}]$. If $\|W_n^\top W_n\|_{\mathrm{F}} > 1$, $\mu_n = 2/\|W_n^\top W_n\|_{\mathrm{F}} \in (0, 2/\|W_n^\top W_n\|_{\mathrm{F}}]$. Accordingly, $\mu_n$ and $\lambda_n$ used in Proposed $(C)$ satisfy steps 4 and 11 in Algorithm \ref{algo:2}.} To see how the choice of parameters $\mu_n$ and $\lambda_n$ $(n\in \mathbb{N})$ affects the convergence rate of Algorithm \ref{algo:2}, Proposed $(C,c)$ with $\mu_n = \lambda_n = c > 0$ $(n\in \mathbb{N})$ was also tested. In Proposed $(C)$, $\mu_n$ and $\lambda_n$ must be computed at each iteration $n$ while in Proposed $(C,c)$ they do because they do not depend on $n$. The  MATLAB source code for Proposed $(C)$ is shown in the Appendix. 

Each of matrices $V := [V_{ij}]\in \mathbb{R}_+^{M \times N}$ in the experiments was generated randomly using the Mersenne Twister Algorithm and the \texttt{rand} function in MATLAB with the rate
\begin{align*}
r := \frac{\# \left\{ V_{ij} \neq 0 \colon 1\leq i \leq M, 1 \leq j \leq N \right\}}{MN},
\end{align*} 
where $\# A$ stands for the number of elements in set $A$ (i.e., all elements of $V$ are nonzero when $r = 1$ and $V$ is a sparse matrix when $r$ is small). One hundred samplings, each starting from a different randomly chosen initial point, were performed until one of the following stopping conditions \cite[{\tt MaxIter}, {\tt TolFun}, {\tt TolX}]{nmf} was satisfied.
\begin{align*}
&n \geq n_{\max} := 10^3, \text{ }
\frac{f\left(W_n, H_n \right) - f \left(W_{n-1}, H_{n-1} \right)}{\max\{1, f \left(W_{n-1}, H_{n-1} \right)\}}
\leq f_{\mathrm{tol}} := 10^{-4},\\
&\max_{\substack{1\leq i \leq M \\ 1 \leq a \leq R}} \frac{\left|\left(W_n - W_{n-1}\right)_{ia}\right|}{\left(W_{n-1} \right)_{ia}} 
\leq x_{\mathrm{tol}} := 10^{-4}, \text{ }
\max_{\substack{1 \leq b \leq R \\ 1 \leq j \leq N}} \frac{\left|\left(H_n - H_{n-1}\right)_{bj}\right|}{\left(H_{n-1} \right)_{bj}} 
\leq x_{\mathrm{tol}} := 10^{-4},
\end{align*} 
where $X_{ij}$ stands for the element in the $i$th row and the $j$th column of a matrix $X := [X_{ij}]$. Four performance measures (metrics) were used:
\begin{itemize}
\item bstT: best computational time [s] 
\item avgT: average computational time [s]
\item bstF: best value of $F$
\item avgF: average value of $F$
\end{itemize}

\subsection{Case where $r = 1$}\label{subsec:4.1}
First, let us consider Problem \eqref{nmf} when $V \in \mathbb{R}_+^{50 \times 25}$ with $r:= 1$ (i.e., 
all elements of $V$ are nonzero) and $R := 5$.

\begin{table}[htb]
\begin{center}
\caption{bstT, avgT, bstF, and avgF for MULT, ALS, Proposed ($C$), and Proposed ($C,c$) algorithms when $V \in \mathbb{R}_+^{50 \times 25}$ with $r:= 1$ and $R := 5$}
\begin{tabular}{l|llll}
\hline
& bstT & avgT & bstF & avgF \\
\hline
MULT & 0.00322128 & 0.00446012 & 0.22998610 & 0.23308853 \\
\hline
ALS & 0.00191543 & 0.00262061 & 0.22810993 & 0.22987739 \\
\hline
Proposed ($0.25$) & 0.00181731 & 0.00281272 & 0.22889555 & 0.23144783 \\
Proposed ($0.50$) & 0.00204749 & 0.00306364 & 0.22904580 & 0.23270967 \\
Proposed ($0.75$) & 0.00313090 & 0.00428935 & 0.23126551 & 0.23559425 \\
\hline
Proposed ($0.25,1$) & 0.00017372 & 0.00038174 & 0.68140744 & 6.45436416 \\
Proposed ($0.50, 1$) & 0.00017256 & 0.00023527 & 0.44183387 & 0.67815429 \\
Proposed ($0.75, 1$) & 0.00020537 & 0.00026407 & 0.33173419 & 0.40918466 \\
\hline
Proposed ($0.25, 0.1$) & 0.00017419 & 0.00093492 & 0.22874089 & 0.30344491 \\
Proposed ($0.50, 0.1$) & 0.00020876 & 0.00218620 & 0.23028413 & 0.24759950 \\
Proposed ($0.75, 0.1$) & 0.00267352 & 0.00337411 & 0.23210681 & 0.23639574 \\
\hline
Proposed ($0.25, 0.01$) & 0.00439306 & 0.00619254 & 0.24045680 & 0.24615864 \\
Proposed ($0.50, 0.01$) & 0.00328624 & 0.00653555 & 0.24445886 & 0.25278974 \\
Proposed ($0.75, 0.01$) & 0.00359608 & 0.00451201 & 0.25839119 & 0.27007882 \\
\hline
\end{tabular}
\label{table:1}
\end{center}
\end{table}

As shown in Table \ref{table:1}, ALS converged faster than MULT, and hence, ALS performed better than MULT. Moreover, Proposed $(0.25)$ converged faster than Proposed $(C)$ $(C = 0.50, 0.75)$. A small parameter value such as $C = 0.25$ (i.e., large coefficients of nonexpansive mappings $T_n$ and $S_n$ in steps 9 and 16 in Algorithm \ref{algo:2}) apparently affects the convergence speed. In particular, Proposed $(0.25)$ (bstT = 0.00181731, avgT = 0.00281272, bstF = 0.22889555, avgF = 0.23144783) had efficiency equivalent to that of ALS (bstT = 0.00191543, avgT = 0.00262061, bstF = 0.22810993, avgF = 0.22987739). Meanwhile, Proposed $(C, 1)$ $(C = 0.25, 0.50, 0.75)$ performed badly. This is because $\mu_n := 1 > 2/\|W_n^\top W_n\|_{\mathrm{F}}$ and $\lambda_n := 1 > 2/\|H_{n+1}^\top H_{n+1}\|_{\mathrm{F}}$ $(n\in \mathbb{N})$ were satisfied, i.e., $T_n$ and $S_n$ in Proposed $(C, 1)$ $(C = 0.25, 0.50, 0.75)$ did not satisfy the nonexpansivity condition. Proposed $(C, c)$ $(C = 0.25, 0.50, 0.75$, $c = 0.1, 0.01)$ performed better than Proposed $(C, 1)$ $(C = 0.25, 0.50, 0.75)$. This implies that a small constant $c$ should be chosen so that $\mu_n := c \leq 2/\|W_n^\top W_n\|_{\mathrm{F}}$ and $\lambda_n := c \leq 2/\|H_{n+1}^\top H_{n+1}\|_{\mathrm{F}}$ $(n\in \mathbb{N})$ are satisfied.

\begin{table}[htb]
\begin{center}
\caption{bstT, avgT, bstF, and avgF for MULT, ALS, Proposed ($C$), and Proposed ($C,c$) algorithms when $V \in \mathbb{R}_+^{100 \times 50}$ with $r:= 1$ and $R := 10$}
\begin{tabular}{l|llll}
\hline
& bstT & avgT & bstF & avgF \\
\hline
MULT & 0.00805119 & 0.00942807 & 0.23493777 & 0.23688365 \\
\hline
ALS & 0.00289740 & 0.00518824 & 0.23454189 & 0.23729251 \\
\hline
Proposed (0.25) & 0.00541920 & 0.00702853 & 0.23331829 & 0.23519642 \\
Proposed (0.50) & 0.00685725 & 0.00847596 & 0.23395753 & 0.23677681 \\
Proposed (0.75) & 0.01014336 & 0.01212223 & 0.23706601 & 0.23995684 \\
\hline
Proposed (0.25, 1) & 0.00028569 & 0.00042164 & 0.59627613 & 1.02473499 \\
Proposed (0.50, 1) & 0.00028367 & 0.00037615 & 0.44638468 & 0.55394564 \\
Proposed (0.75, 1) & 0.00036966 & 0.00054351 & 0.34954833 & 0.49603678 \\
\hline
Proposed (0.25, 0.1) & 0.00028824 & 0.00039870 & 0.43310099 & 0.47716912 \\
Proposed (0.50, 0.1) & 0.00037289 & 0.00048605 & 0.28937148 & 0.31896510 \\
Proposed (0.75, 0.1) & 0.00660396 & 0.00752956 & 0.23499926 & 0.23671299 \\
\hline
Proposed (0.25, 0.01) & 0.01578557 & 0.01804486 & 0.24512709 & 0.24943914 \\
Proposed (0.50, 0.01) & 0.01434792 & 0.01801145 & 0.25084821 & 0.25655329 \\
Proposed (0.75, 0.01) & 0.00748992 & 0.00850517 & 0.27290841 & 0.27580332 \\
\hline
\end{tabular}
\label{table:2}
\end{center}
\end{table}

Next, let us consider Problem \eqref{nmf} when $V \in \mathbb{R}_+^{100 \times 50}$ with $r:= 1$ and $R := 10$. From Table \ref{table:2}, MULT, ALS, and Proposed $(C)$ $(C = 0.25, 0.50, 0.75)$ converged to almost the same value of $F$. Although Proposed (0.25) (bstT = 0.00541920, avgT = 0.00702853) converged faster than MULT (bstT = 0.00805119, avgT = 0.00942807), ALS (bstT = 0.00289740, avgT = 0.00518824) converged faster than Proposed (0.25). Meanwhile, Proposed $(C)$ $(C = 0.25, 0.50)$ (e.g., Proposed (0.25) had bstF = 0.23331829 and avgF = 0.23519642) better optimized function $F$ in Problem \eqref{nmf} than MULT (bstF = 0.23493777, avgF = 0.23688365) and ALS (bstF= 0.23454189, avgF = 0.23729251). Although the efficiency of Proposed $(0.75, 0.1)$ was almost the same as that of Proposed $(C)$ $(C = 0.25, 0.50)$, it is likely that Proposed $(C,0.01)$ performed better than Proposed $(C,c)$ $(c= 1,0.1)$, as also seen in Table \ref{table:1}.

\begin{table}[htb]
\begin{center}
\caption{bstT, avgT, bstF, and avgF for MULT, ALS, Proposed ($C$), and Proposed ($C,c$) algorithms when $V \in \mathbb{R}_+^{200 \times 100}$ with $r:= 1$ and $R := 20$}
\begin{tabular}{l|llll}
\hline
& bstT & avgT & bstF & avgF \\
\hline
MULT & 0.02400269 & 0.02734322 & 0.23650149 & 0.23768112 \\
\hline
ALS & 0.00697689 & 0.01822636 & 0.24222807 & 0.24576315 \\
\hline
Proposed (0.25) & 0.02531202 & 0.02815573 & 0.23560407 & 0.23648003 \\
Proposed (0.50) & 0.03049898 & 0.03466635 & 0.23726896 & 0.23819649 \\
Proposed (0.75) & 0.04484526 & 0.04831823 & 0.24092255 & 0.24219488 \\
\hline
Proposed (0.25, 1) & 0.00070745 & 0.00085857 & 1.73067788 & 1.86318161 \\
Proposed (0.50, 1) & 0.00099585 & 0.00112558 & 0.48901480 & 1.41063377 \\
Proposed (0.75, 1) & 0.00153680 & 0.00171749 & 0.44464377 & 0.68612137 \\
\hline
Proposed (0.25, 0.1) & 0.00070265 & 0.00083291 & 0.34391009 & 0.35276281 \\
Proposed (0.50, 0.1) & 0.00101222 & 0.00118885 & 0.33204987 & 0.34571235 \\
Proposed (0.75, 0.1) & 0.00217970 & 0.02127195 & 0.23459599 & 0.24120886 \\
\hline
Proposed (0.25, 0.01) & 0.00099478 & 0.00114034 & 0.67804880 & 0.69957591 \\
Proposed (0.50, 0.01) & 0.05605289 & 0.06314479 & 0.24633999 & 0.24913422 \\
Proposed (0.75, 0.01) & 0.05598053 & 0.06748969 & 0.25404490 & 0.25782798 \\
\hline
\end{tabular}
\label{table:3}
\end{center}
\end{table}

Table \ref{table:3} shows the metrics for Problem \eqref{nmf} when $V \in \mathbb{R}_+^{200 \times 100}$ with $r := 1$ and $R := 20$. Although MULT (bstT = 0.02400269, avgT = 0.02734322) and ALS (bstT = 0.00697689, avgT = 0.01822636) converged faster than Proposed $(C)$ $(C = 0.25, 0.50, 0.75)$, Proposed $(0.25)$ (bstF = 0.23560407, avgF = 0.23648003) optimized function $F$ better than MULT (bstF = 0.23650149, avgF = 0.23768112) and ALS (bstF = 0.24222807, avgF = 0.24576315), as seen in Table \ref{table:2}.

These results show that, when all the elements of $V$ are nonzero and $MN$ is small, the efficiency of Proposed $(C)$ with small constant $C$ is almost the same as that of ALS, which is the best algorithm for NMF. When all the elements of $V$ are nonzero and $MN$ is large, ALS converges faster than the proposed algorithms while Proposed $(C)$ with small constant $C$ optimizes $F$ better than MULT and ALS.

\subsection{Case where $r \leq 0.01$}\label{subsec:4.2}

\begin{table}[htb]
\begin{center}
\caption{bstT, avgT, bstF, and avgF for MULT, ALS, Proposed ($C$), and Proposed ($C,c$) algorithms when $V \in \mathbb{R}_+^{50 \times 25}$ with $r \leq 0.01$ and $R := 5$}
\begin{tabular}{l|llll}
\hline
& bstT & avgT & bstF & avgF \\
\hline
MULT & 0.00166420 & 0.00232842 & 0.02751133 & 0.02937415 \\
\hline
ALS & 0.00165955 & 0.00205977 & 0.02750745 & 0.02863840 \\
\hline
Proposed (0.25) & 0.00070293 & 0.00113707 & 0.02752362 & 0.02909619 \\
Proposed (0.50) & 0.00071851 & 0.00118315 & 0.02755326 & 0.02954198 \\
Proposed (0.75) & 0.00117840 & 0.00193609 & 0.02766379 & 0.02978319 \\
\hline
Proposed (0.25, 1) & 0.00059409 & 0.00091400 & 0.02755416 & 0.02904506 \\
Proposed (0.50, 1) & 0.00076439 & 0.00125031 & 0.02763945 & 0.02929433 \\
Proposed (0.75, 1) & 0.00119457 & 0.00168017 & 0.02777880 & 0.02990826 \\
\hline
Proposed (0.25, 0.1) & 0.00026702 & 0.00041036 & 0.05665250 & 0.05756447 \\
Proposed (0.50, 0.1) & 0.00041605 & 0.00058405 & 0.05610456 & 0.05685127 \\
Proposed (0.75, 0.1) & 0.00071190 & 0.00092664 & 0.05644580 & 0.05715633 \\
\hline
Proposed (0.25, 0.01) & 0.00104115 & 0.00141259 & 0.05746258 & 0.05793165 \\
Proposed (0.50, 0.01) & 0.00153107 & 0.00202124 & 0.05781522 & 0.05844957 \\
Proposed (0.75, 0.01) & 0.00281485 & 0.00357462 & 0.05897798 & 0.06012710 \\
\hline
\end{tabular}
\label{table:4}
\end{center}
\end{table}

\begin{table}[htb]
\begin{center}
\caption{bstT, avgT, bstF, and avgF for MULT, ALS, Proposed ($C$), and Proposed ($C,c$) algorithms when $V \in \mathbb{R}_+^{100 \times 50}$ with $r \leq 0.01$ and $R := 10$}
\begin{tabular}{l|llll}
\hline
& bstT & avgT & bstF & avgF \\
\hline
MULT & 0.00266517 & 0.00381133 & 0.03161738 & 0.03319867 \\
\hline
ALS & 0.00255122 & 0.00318825 & 0.03159433 & 0.03270456 \\
\hline
Proposed (0.25) & 0.00153508 & 0.00250549 & 0.03169617 & 0.03349920 \\
Proposed (0.50) & 0.00206164 & 0.00328698 & 0.03187789 & 0.03376532 \\
Proposed (0.75) & 0.00347115 & 0.00526441 & 0.03229713 & 0.03468006 \\
\hline
Proposed (0.25, 1) & 0.00131476 & 0.00181448 & 0.03161219 & 0.03289926 \\
Proposed (0.50, 1) & 0.00172756 & 0.00238345 & 0.03176905 & 0.03322947 \\
Proposed (0.75, 1) & 0.00306115 & 0.00414748 & 0.03183166 & 0.03374404 \\
\hline
Proposed (0.25, 0.1) & 0.00053445 & 0.00074888 & 0.05831575 & 0.05888618 \\
Proposed (0.50, 0.1) & 0.00094323 & 0.00108989 & 0.05852361 & 0.05877634 \\
Proposed (0.75, 0.1) & 0.00177545 & 0.00198634 & 0.05875030 & 0.05901237 \\
\hline
Proposed (0.25, 0.01) & 0.00069307 & 0.00083472 & 0.05905958 & 0.05920412 \\
Proposed (0.50, 0.01) & 0.00119537 & 0.00150471 & 0.05920358 & 0.05941218 \\
Proposed (0.75, 0.01) & 0.00323848 & 0.00400472 & 0.05990993 & 0.06024958 \\
\hline
\end{tabular}
\label{table:5}
\end{center}
\end{table}

Table \ref{table:4} shows the metrics for Problem \eqref{nmf} when $R := 5$ and $V \in \mathbb{R}_+^{50 \times 25}$ with $r \leq 0.01$, i.e., the number of elements in $V$ is $1250$ and $\# \{V_{ij} \neq 0 \colon 1\leq i \leq 50, 1 \leq j \leq 25 \} \leq 12.5$.

The results in Table \ref{table:4} show that Proposed ($C$) $(C = 0.25, 0.50, 0.75)$ (e.g., Proposed (0.25) had bstT = 0.00070293 and avgT = 0.00113707) converges faster than MULT (bstT = 0.00166420, avgT = 0.00232842) and ALS (bstT = 0.00165955, avgT = 0.00205977) and that bstF and the avgF for Proposed (0.25) (bstF = 0.02752362, avgF = 0.02909619) was almost the same as the bstF and avgF for ALS (bstF = 0.02750745, avgF = 0.02863840). 
Moreover, Proposed (0.25,1) converged the fastest of all the algorithms used in this experiment,
and 
the proposed algorithms with $C = 0.25$ performed better than the ones with $C = 0.50, 0.75$. 
This implies that the larger the weighted parameters of nonexpansive mappings $T_n$ and $S_n$, the shorter the computational time. In contrast to the results in subsection \ref{subsec:4.1}, Proposed ($C,1$) $(C = 0.25, 0.50, 0.75)$ performed better than Proposed ($C,c$) $(C = 0.25, 0.50, 0.75, c = 0.1, 0.01)$. When $V$ is sparse, $\mu_n = \lambda_n = c$ $(c=1,0.1, 0.01)$ satisfies $\mu_n \leq 2/\|W_n^\top W_n \|_{\mathrm{F}}$ and $\lambda_n \leq 2/\|H_{n+1}^\top H_{n+1} \|_{\mathrm{F}}$. However, if $c$ is small and if $W_n, H_n \approx O$, then $T_n (H_n) = P_{\mathbb{R}_+^{R \times N}}[H_n - \mu_n W_n^\top (W_n H_n - V)] \approx P_{\mathbb{R}_+^{R \times N}}(H_n) = H_n$ and $S_n (W_n) = P_{\mathbb{R}_+^{M \times R}} [W_n - \lambda_n (W_n H_{n+1} - V) H_{n+1}^\top] \approx P_{\mathbb{R}_+^{M \times R}} (W_n) = W_n$ hold in the sense of the Frobenius norm, which implies that the gradient of $F$ has little or no effect. This is one reason the proposed algorithms with small constant $c$ could not optimize $F$.

Next, let us consider Problem \eqref{nmf} when $R := 10$ and $V \in \mathbb{R}_+^{100 \times 50}$ with $r \leq 0.01$, i.e., the number of the elements in $V$ is $5000$ and $\# \{V_{ij} \neq 0 \colon 1\leq i \leq 100, 1 \leq j \leq 50 \} \leq 50$. Table \ref{table:5} shows that Proposed (0.25) and Proposed (0.25, 1) converged faster than MULT and ALS and that Proposed ($C,1$) $(C = 0.25, 0.50, 0.75)$ performed better than Proposed ($C,c$) $(C = 0.25, 0.50, 0.75, c = 0.1, 0.01)$, as seen in Table \ref{table:4}. 

\begin{table}[htb]
\begin{center}
\caption{bstT, avgT, bstF, and avgF for MULT, ALS, Proposed ($C$), and Proposed ($C,c$) algorithms when $V \in \mathbb{R}_+^{200 \times 100}$ with $r \leq 0.01$ and $R := 20$}
\begin{tabular}{l|llll}
\hline
& bstT & avgT & bstF & avgF \\
\hline
MULT & 0.00500478 & 0.00731791 & 0.03988082 & 0.04069823 \\
\hline
ALS & 0.00696703 & 0.00948289 & 0.03964874 & 0.04020298 \\
\hline
Proposed (0.25) & 0.00760513 & 0.00995760 & 0.04036087 & 0.04127604 \\
Proposed (0.50) & 0.01015383 & 0.01306092 & 0.04049975 & 0.04173487 \\
Proposed (0.75) & 0.01449030 & 0.01944444 & 0.04090514 & 0.04305908 \\
\hline
Proposed (0.25, 1) & 0.00403276 & 0.00514699 & 0.03972592 & 0.04041520 \\
Proposed (0.50, 1) & 0.00560411 & 0.00750421 & 0.03988278 & 0.04061267 \\
Proposed (0.75, 1) & 0.01142270 & 0.01364317 & 0.04020438 & 0.04115410 \\
\hline
Proposed (0.25, 0.1) & 0.00155693 & 0.00203610 & 0.05974672 & 0.05995876 \\
Proposed (0.50, 0.1) & 0.00300555 & 0.00327410 & 0.05961122 & 0.05967715 \\
Proposed (0.75, 0.1) & 0.00539418 & 0.00590893 & 0.05980458 & 0.05990551 \\
\hline
Proposed (0.25, 0.01) & 0.00182208 & 0.00199416 & 0.06000842 & 0.06014881 \\
Proposed (0.50, 0.01) & 0.00326519 & 0.00357171 & 0.05969398 & 0.05981723 \\
Proposed (0.75, 0.01) & 0.00861478 & 0.00968720 & 0.06054412 & 0.06070213 \\
\hline
\end{tabular}
\label{table:6}
\end{center}
\end{table}

Table \ref{table:6} shows the metrics for Problem \eqref{nmf} when $R := 20$ and $V\in \mathbb{R}_+^{200 \times 100}$ with $r \leq 0.01$, i.e., the number of elements in $V$ is 20000 and $\# \{V_{ij} \neq 0 \colon 1 \leq i \leq 200, 1 \leq j \leq 100 \} \leq 200$. In contrast to Table \ref{table:5}, MULT converged faster than ALS. Moreover, Proposed (0.25, 1) (bstT = 0.00403276, avgT = 0.00514699) converged faster than MULT (bstT = 0.00500478, avgT = 0.00731791), and the value of $F$ for Proposed (0.25, 1) (bstF = 0.03972592, avgF = 0.04041520) were almost the same as the one for MULT (bstF = 0.03988082, avgF = 0.04069823).

From these results, we conclude that, when $V$ is sparse, Proposed $(C)$ and Proposed ($C,c$) with small constant $C$ and large constant $c$ perform better than MULT and ALS and that the proposed algorithms with large constant $c$ perform better than the ones with small constant $c$.

\subsection{Case where $r \leq 0.001$}\label{subsec:4.3}  
Finally, we apply MULT, ALS, and Algorithm \ref{algo:2} to Problem \eqref{nmf} when $R := 20$ and $V \in \mathbb{R}_+^{200 \times 100}$ with $r \leq 0.001$, i.e., the number of the elements in $V$ is $20000$ and $\# \{V_{ij} \neq 0 \colon 1\leq i \leq 200, 1 \leq j \leq 100 \} \leq 20$.

\begin{table}[htb]
\begin{center}
\caption{bstT, avgT, bstF, and avgF for MULT, ALS, Proposed ($C$), and Proposed ($C,c$) algorithms when $V \in \mathbb{R}_+^{200 \times 100}$ with $r \leq 0.001$ and $R := 20$}
\begin{tabular}{l|llll}
\hline
& bstT & avgT & bstF & avgF \\
\hline
MULT & 0.00500909 & 0.00717753 & 0.00055390 & 0.00287419 \\
\hline
ALS & 0.00374294 & 0.00468982 & 0.00000000 & 5.04918745 \\
\hline
Proposed (0.25) & 0.00732760 & 0.01085649 & 0.00176924 & 0.00310601 \\
Proposed (0.50) & 0.00836821 & 0.01321850 & 0.00195163 & 0.00375576 \\
Proposed (0.75) & 0.01350137 & 0.01875023 & 0.00261926 & 0.00518008 \\
\hline
Proposed (0.25, 1) & 0.00713262 & 0.00976884 & 0.00169725 & 0.00299880 \\
Proposed (0.50, 1) & 0.00940968 & 0.01283440 & 0.00202390 & 0.00353762 \\
Proposed (0.75, 1) & 0.01522651 & 0.02016161 & 0.00246735 & 0.00515722 \\
\hline
Proposed (0.25, 0.1) & 0.00182897 & 0.00202514 & 0.01725790 & 0.01727294 \\
Proposed (0.50, 0.1) & 0.00300825 & 0.00323018 & 0.01724517 & 0.01726637 \\
Proposed (0.75, 0.1) & 0.00622978 & 0.00706693 & 0.01743940 & 0.01749223 \\
\hline
Proposed (0.25, 0.01) & 0.00209641 & 0.00230292 & 0.01725625 & 0.01727153 \\
Proposed (0.50, 0.01) & 0.00347790 & 0.00374012 & 0.01725655 & 0.01731079 \\
Proposed (0.75, 0.01) & 0.01075377 & 0.01297691 & 0.01830596 & 0.01847198 \\
\hline
\end{tabular}
\label{table:7}
\end{center}
\end{table}

As shown in Table \ref{table:7}, the convergence of ALS depended on the initial points. Sometimes ALS found the global minimizer of $F$ (bstF = 0.00000000), and sometimes ALS converged to a point ranked lower than $R$, which is not optimal (avgF = 5.04918745). Meanwhile, MULT, Proposed (0.25), and Proposed (0.25,1) were robust for Problem \eqref{nmf} with the sparse matrix $V$, as seen in Tables \ref{table:4}--\ref{table:6}. Table \ref{table:7} also shows that MULT converged faster than Proposed ($C,c$) $(C = 0.25, 0.50, 0.75, c = 1, 0.1, 0.01)$. This is because the setting of $C = 0.25$ and $c =1$ is insufficient for Problem \eqref{nmf} when $r$ is too small. Therefore, on the basis of the discussion in subsection \ref{subsec:4.2}, we checked whether Proposed ($C,c$) with $C < 0.25$ and $c > 1$ converges faster than MULT. We used $C = 0.20$ and $c = 2$ for example.

\begin{table}[htb]
\begin{center}
\caption{bstT, avgT, bstF, and avgF for MULT, ALS, and Proposed ($0.20, 2$) algorithms when $V \in \mathbb{R}_+^{200 \times 100}$ with $r \leq 0.001$ and $R := 20$}
\begin{tabular}{l|llll}
\hline
& bstT & avgT & bstF & avgF \\
\hline
MULT & 0.00500909 & 0.00717753 & 0.00055390 & 0.00287419 \\
\hline
ALS & 0.00374294 & 0.00468982 & 0.00000000 & 5.04918745 \\
\hline
Proposed (0.20, 2) & 0.00363362 & 0.00549641 & 0.00062403 & 0.00253986\\
\hline
\end{tabular}
\label{table:8}
\end{center}
\end{table}

As shown in Table \ref{table:8}, Proposed (0.20, 2) (bstT = 0.00363362, avgT = 0.00549641) converged faster than MULT (bstT = 0.00500909, avgT = 0.00717753). Therefore, we can conclude that Proposed ($C$) and Proposed ($C,c$) with small constant $C$ and large constant $c$ should be used for Problem \eqref{nmf} when $V$ is sparse, as also stated in subsection \ref{subsec:4.2}.

\section{Conclusion and future work}\label{sec:5}
We considered the nonmonotone variational inequality in NMF and presented a fixed point algorithm for solving the variational inequality. Investigation of the convergence property of the proposed algorithm showed that, under certain assumptions, any accumulation point of the sequence generated by the proposed algorithm belongs to the solution set of the variational inequality.

We also numerically compared the proposed algorithm with the {\tt 'mult'} and {\tt 'als'} algorithms for NMF optimization problems in certain situations. Numerical examples demonstrated that, when all the elements of a data matrix are nonzero, the efficiency of the proposed algorithm is almost the same as that of the {\tt 'als'} algorithm and that, when a data matrix is sparse, the proposed algorithm converges faster than the {\tt 'mult'} and {\tt 'als'} algorithms. Moreover, the {\tt 'als'} algorithm may not converge to an optimal solution. The numerical examples also suggested parameter designs that ensure fast convergence of the proposed algorithm. In particular, the weighted parameters of the nonexpansive mappings used in the proposed algorithm should be large to solve NMF problems effectively. 

Our main objective was to devise a fixed point algorithm for NMF based on the Krasnosel'ski\u\i-Mann fixed point algorithm. Another particularly interesting problem is determining whether there are fixed point algorithms for NMF based on the Halpern fixed point algorithm and the hybrid method. We will thus consider developing fixed point algorithms for NMF based on the Halpern fixed point algorithm and the hybrid method and evaluating their capabilities.

\section*{Appendix}\label{appendix}
\lstlistingname~\ref{lst:nmfmt} is the MATLAB source code for Proposed $(C)$ (Algorithm \ref{algo:2}) used in the experiments (see \cite{nmf} for the methods for generating $(W^\star, H^\star)$ satisfying (i) and (ii) in Problem \eqref{nmf}). The function in \lstlistingname~\ref{lst:nmfmt} requires five parameters: V, W, H, C, and MAX\_ITER. Parameter V is a factorized nonnegative matrix, and W and H stand for the initial nonnegative matrices in Algorithm \ref{algo:2}. Parameter C is a ratio of convex combination, which implies that C is equal to $\alpha_n = \lambda_n \in (0,1)$ $(n\in \mathbb{N})$ in Algorithm \ref{algo:2}. Parameter MAX\_ITER is the maximum number of iterations needed to compute Algorithm \ref{algo:2}. Other stopping conditions can be included in \textcolor{comment}{\%\{Stopping Criteria\}\%} in \lstlistingname~\ref{lst:nmfmt}.

\begin{lstlisting}[caption=MATLAB source code for Proposed ($C$),label=lst:nmfmt]
function [W, H, D] = nmfmt(V, W, H, C, MAX_ITER)

NM = numel(V);
D = Inf;
for it = 1:MAX_ITER
    % Update for H
    st = 2. / max(1., norm(W.' * W, 'fro'));
    T = max(0, H - st * W.' * (W * H - V));
    H = C * H + (1. - C) * T;
    % Update for W
    st = 2. / max(1., norm(H.' * H, 'fro'));
    T = max(0, W - st * (W * H - V) * H.');
    W = C * W + (1. - C) * T;
    % Evaluate current solution
    DAX = V - W * H;
    D = sqrt(sum(sum(DAX.^2)) / NM);
    if %{ Stopping Criteria }%
        break;
    end
end
% Normalization of W and H
end
\end{lstlisting}

The function in \lstlistingname~\ref{lst:nmfmt} returns three values: W, H, and D.
(W, H) is a solution to Problem \eqref{nmf}, and 
D stands for the value of $F$ in Problem \eqref{nmf} at (W, H).\\

\textbf{Acknowledgments} 
We thank Kazuhiro Hishinuma for his input on the numerical examples.

\end{document}